\documentclass[11pt,a4paper]{article}
\date{}
\makeatletter
\@addtoreset{equation}{section}
\makeatother

\usepackage[colorlinks,
            linkcolor=blue,
            anchorcolor=blue,
            citecolor=blue]{hyperref}
\usepackage{hyperref}
\hypersetup{colorlinks=true,
            linkcolor=blue,
            anchorcolor=blue,
            citecolor=blue}
\usepackage{hyperref,amsmath,amssymb,amscd}
\usepackage{amssymb,amsmath,enumerate}
\usepackage{indentfirst}
\usepackage{url}
\usepackage[capitalise]{cleveref}
\usepackage{amsthm}
\usepackage{amsmath}
\AtBeginDocument{\hypersetup{pdfborder={0 0 0.01}}}
\newtheorem{theorema}{Theorem}[section]

\newtheorem{thm}{Theorem}[section]
\newtheorem{cor}[thm]{Corollary}
\newtheorem{lem}[thm]{Lemma}

\theoremstyle{definition}

\newtheorem{rem}[thm]{\bf Remark}
\allowdisplaybreaks[4]

\usepackage{enumitem}
\usepackage{extarrows}
\usepackage{amsthm}
\usepackage{setspace}
\hypersetup{urlcolor=blue, citecolor=red}

\begin{document}

\title{\bf
Remainder terms of $L^p$-Hardy inequalities with magnetic fields: the case $1<p<2$\footnotetext{This work was supported by National Natural Science Foundation of China (No. 12371120) and Chongqing Graduate Student Research Innovation Project (No. CYB240097).}}\author{{Xiao-Ping Chen,\ \ Chun-Lei Tang\footnote{Corresponding author.\newline
\indent\,\,\, \emph{E-mail address:} tangcl@swu.edu.cn (C.-L. Tang); xpchen\_maths@163.com (X.-P. Chen).}}\\
{\small \emph{School of Mathematics and Statistics, Southwest University,  Chongqing {\rm400715},}}\\
{\small \emph{People's Republic of China}}}
\maketitle
\baselineskip 17pt

\noindent {\bf Abstract}:\ This paper focuses on remainder estimates of the magnetic $L^p$-Hardy inequalities for $1<p<2$. \emph{Firstly}, we establish a family of remainder terms involving magnetic gradients of the magnetic $L^p$-Hardy inequalities, which are also new even for the classical $L^p$-Hardy inequalities. \emph{Secondly}, we study another family of remainder terms involving logarithmic terms of the magnetic $L^p$-Hardy inequalities. \emph{Lastly}, as a byproduct, we further obtain remainder terms of some other $L^p$-Hardy-type inequalities by using similar proof of our main results.

Furthermore, this paper answers the open question proposed by Cazacu \emph{et al.} in [Nonlinearity 37:035004, 2024] and can be viewed as a supplementary work of it.

\vspace{0.25em}


\vspace{0.25em}


\section{Introduction and main results}\label{sect-1}

\noindent The purpose of this paper is to study remainder estimates of $L^p$-Hardy inequalities with magnetic fields for $1<p<2$. The novelty of this paper is stated as follows.
\begin{enumerate}
[itemsep=0pt, topsep=2pt, parsep=0pt]

\item[$(a)$]
By adding remainder terms involving magnetic gradients or logarithmic terms, we improve the magnetic $L^p$-Hardy inequalities for $1<p<2$.

\item[$(b)$]
Compared with \cite{Cazacu24} where the authors study remainder terms of the magnetic $L^p$-Hardy inequalities for $p\ge2$, this paper fills it for the case $1<p<2$.

\item[$(c)$]
This paper replies the open question proposed in \cite[Remark 1.1]{Cazacu24}, and can be viewed as a supplementary work of it.

\item[$(d)$]
Using similar proof of our main results, we further establish remainder terms of some other $L^p$-Hardy-type inequalities.

\end{enumerate}

We recall some enlightening results of Hardy-type inequalities in Section \ref{sect-1.1}. In Section \ref{sect-1.2}, we describe main results and some comments. Section \ref{sect-1.3} contains the outline of the rest of this paper.

\subsection{Overview and motivation}\label{sect-1.1}

\subsubsection{The classical $L^p$-Hardy inequalities}
\label{sect-1.1.1}

\noindent The classical $L^p$-Hardy inequalities state as follows: for $1<p<N$,
\begin{equation}\label{1.2}
\int_{\mathbb{R}^N}|\nabla u|^p\mathrm{d}x
\ge\left(\frac{N-p}{p}\right)^p
\int_{\mathbb{R}^N}\frac{|u|^p}{|x|^p}\mathrm{d}x,
\ \
\forall
u\in C_0^\infty(\mathbb{R}^N\setminus\{\mathbf{0}\}),
\end{equation}
see \cite{Hardy88} for details. Moreover, the constant $\left(\frac{N-p}{p}\right)^p$ is optimal, but there exists no nontrivial function such that the equality of \eqref{1.2} occurs. In the study of elliptic and parabolic equations, the classical $L^p$-Hardy inequalities \eqref{1.2} and their improved versions serve as  indispensable tools (see \cite{Brezis97,Fleckinger99} \emph{etc.}).

In the present paper, we are interested in remainder estimates of $L^p$-Hardy-type inequalities. It is worth pointing out that the analysis of remainder estimates has been given a lot of attention after the open question proposed by Brezis and Lieb in \cite{Brezis85}. Next, we describe the research status of remainder estimates for $L^p$-Hardy-type inequalities.

For the special case $p=N=2$, under some extra assumptions on $u$, there still exist some nontrivial Hardy-type inequalities for the operator $-\Delta$ when the weight $\frac{1}{|x|^2}$ instead of a weaker one involving an additional logarithmic term. Namely,
\begin{equation*}
\int_{\mathbb{R}^2}|\nabla u|^2\mathrm{d}x
\ge C\int_{\mathbb{R}^2}
\frac{|u|^2}{|x|^2(1+|\ln |x||^2)}\mathrm{d}x,
\ \
\mathrm{if}
\
\int_{\{|x|=1\}}u(x)\mathrm{d}x=0.
\end{equation*}
See for example, \cite{Laptev98,Solomyak94}. For $2\le p<N$, it has been established in \cite{Sano18} that
$$
\int_{\mathbb{R}^N}
\left|\frac{x}{|x|}\cdot\nabla u\right|^p\mathrm{d}x
-\left(\frac{N-p}{p}\right)^p
\int_{\mathbb{R}^N}\frac{|u|^p}{|x|^p}\mathrm{d}x
\ge C\sup\limits_{R>0}\int_{\mathbb{R}^N}
\frac{|u(x)-v(x)|^p}
{|x|^p\left|\ln\frac{R}{|x|}\right|^p}
\mathrm{d}x,
$$
for some $C=C(p,N)>0$, where $v(x)=|x|^{\frac{p-N}{p}}
R^{\frac{N-p}{p}}u\left(R\frac{x}{|x|}\right)
$. To obtain more related results about remainder estimates of $L^p$-Hardy-type inequalities (including bounded domains), we refer to \cite{Brezis97,Cianchi08,Dyda24,Filippas02,Ruzhansky18} and the reference therein.

\subsubsection{The magnetic $L^p$-Hardy inequalities}
\label{sect-1.1.2}

\noindent We say that the \emph{magnetic tensor} $\mathbf{B}:\mathbb{R}^N\to\mathbb{R}^{N\times N}$ ($2$-form) is exact if there exists a \emph{magnetic potential} $\mathbf{A}:\mathbb{R}^N\to\mathbb{R}^N$ (smooth $1$-form) satisfying $\mathrm{d}\mathbf{A}=\mathbf{B}$. In particular, $\mathbf{B}$ is closed, that is, $\mathrm{d}\mathbf{B}=0$, where $``\mathrm{d}"$ is the exterior derivative. From the Poincar\'{e} lemma (see \cite[Corollary 18]{Spivak79}), $\mathbf{B}$ is closed if and only if $\mathbf{B}$ is exact.

We define the magnetic $p$-Laplacian
\begin{equation*}
\Delta_{\mathbf{A},p}u:
=\mathrm{div}\!_\mathbf{A}
\left(\left|\nabla\!_\mathbf{A}u\right|^{p-2}
\nabla\!_\mathbf{A}u\right)
\end{equation*}
on $C_0^\infty(\mathbb{R}^N)$, where $\mathrm{div}\!_\mathbf{A}f
:=\mathrm{div}f+i\mathbf{A}\cdot f$ is the magnetic divergence and $\nabla\!_\mathbf{A}u:=\nabla u+i\mathbf{A}(x)u$ is the magnetic gradient. Obviously, $\Delta_{0,p}=\Delta_{p}$.

We also denote
\[
\mathcal{D}_{\mathbf{A},p}
:=\overline{C_0^\infty
(\mathbb{R}^N\setminus\{\mathbf{0}\})}
^{\|\cdot\|},
\]
with the norm
\begin{equation*}
\|u\|:=\left(
\int_{\mathbb{R}^N}
\left|\nabla\!_\mathbf{A} u\right|^p\mathrm{d}x
+\int_{\mathbb{R}^N}
\left|u\right|^p
\mathrm{d}x\right)^{\frac{1}{p}}.
\end{equation*}
Furthermore, from \cite[Theorem 7.21]{Lieb97}, the diamagnetic inequality
\begin{equation}\label{1.1}
\left|\nabla\!_\mathbf{A} u\right|
\ge\left|\nabla|u|\right|
\end{equation}
holds for $a.e.$ $x\in\mathbb{R}^N$ and all $u\in\mathcal{D}_{\mathbf{A},p}$. Then, due to \eqref{1.2},  if $1<p<N$,
\begin{equation*}
\int_{\mathbb{R}^N}
\left|\nabla\!_\mathbf{A} u\right|^p\mathrm{d}x
\ge\int_{\mathbb{R}^N}
\left|\nabla|u|\right|^p\mathrm{d}x
\ge\left(\frac{N-p}{p}\right)^p
\int_{\mathbb{R}^N}\frac{|u|^p}{|x|^p}\mathrm{d}x.
\end{equation*}
For simplicity, let us denote
\begin{equation}\label{1.6}
\mathcal{H}_{\mathbf{A},p}(u):
=\int_{\mathbb{R}^N}
\left|\nabla\!_\mathbf{A} u\right|^p\mathrm{d}x
-\left(\frac{N-p}{p}\right)^p
\int_{\mathbb{R}^N}\frac{|u|^p}{|x|^p}\mathrm{d}x
\ge0.
\end{equation}

Now, we present some results about remainder estimates of \eqref{1.6}. When $p=N=2$, Laptev and Weidl in \cite{Laptev98} obtained that, for all smooth $\mathbf{A}:\mathbb{R}^2\to\mathbb{R}^2$ satisfying $\mathbf{B}=\mathrm{curl}\mathbf{A}$,
\[
\int_{\mathbb{R}^2}|\nabla\!_\mathbf{A}u|^2\mathrm{d}x
\ge C(\mathbf{B})
\int_{\mathbb{R}^2}\frac{|u|^2}{1+|x|^2}\mathrm{d}x
\]
provided that
\[
\frac{1}{2\pi}
\int_{\{|x|<r\}}\mathbf{B}\mathrm{d}x
\]
is not an integer. Recently, with no other assumptions on magnetic field $\mathbf{B}\neq0$, it was established by Cazacu and Krej\v{c}i\v{r}\'{i}k in \cite{Cazacu16} that, when $N\ge p=2$, for all smooth $\mathbf{A}$ such that $\mathrm{d}\mathbf{A}=\mathbf{B}$,
\[
\mathcal{H}_{\mathbf{A},2}(u)
\ge C(\mathbf{B},N)\int_{\mathbb{R}^N}
\frac{|u|^2}{1+|x|^2|\ln |x||^2}\mathrm{d}x.
\]
For more related results of $L^p$-Hardy-type inequalities, we refer to \cite{Balinsky04,Ekholm14,Kovarik11} in general magnetic fields, \cite{Alziary03,Fanelli20,Lam23,Lu24} in Aharonov-Bohm magnetic fields (a special case of magnetic fields), \cite{Cassano23} in the Heisenberg groups, \cite{Evans05,Lam23} for Rellich inequalities, \emph{etc}.

More importantly, Cazacu, Krej\v{c}i\v{r}\'{i}k, Lam and Laptev in \cite{Cazacu24} obtained the following results, which improve the results presented in \cite{Cazacu16} from $L^2$ to general $L^p$.

\begin{theorema}
[{\rm{\cite[Theorem 1.3]{Cazacu24}}}]
\label{thm-a}
Assume that $2\le p<N$, $\mathbf{B}$ is nontrivial, smooth and closed. Then for all $\mathbf{A}$ satisfying $\mathrm{d}\mathbf{A}=\mathbf{B}$,
\begin{equation*}
\mathcal{H}_{\mathbf{A},p}(u)
\ge C(\mathbf{B},p,N)\int_{\mathbb{R}^N}
\frac{|u|^p}{|x|^p\left(1+|\ln|x||^p\right)}
\mathrm{d}x,
\ \ \
\forall
u\in\mathcal{D}_{\mathbf{A},p},
\end{equation*}
for some constant $C(\mathbf{B},p,N)>0$.
\end{theorema}

Here we point out that the above result is heavily dependent on the following result, which is another type remainder terms of \eqref{1.6} for $p\ge2$.

\begin{theorema}
[{\rm{\cite[Theorem 1.2]{Cazacu24}}}]
\label{thm-b}
Assume that $2\le p<N$, $\mathbf{B}$ is nontrivial, smooth and closed. Then for all $\mathbf{A}$ satisfying $\mathrm{d}\mathbf{A}=\mathbf{B}$,
\begin{equation*}
\mathcal{H}_{\mathbf{A},p}(u)
\ge c(p)\int_{\mathbb{R}^N}
|x|^{p-N}
\left|\nabla\!_\mathbf{A} \left(u|x|^{\frac{N-p}{p}}\right)\right|^p
\mathrm{d}x,
\ \
\forall
u\in\mathcal{D}_{\mathbf{A},p},
\end{equation*}
where
\[
c(p):=\inf_{s^2+t^2>0}
\frac{\left(t^2+s^2+2s+1\right)^{\frac{p}{2}}-1-ps}
{(t^2+s^2)^{\frac{p}{2}}}
\in(0,1].
\]
\end{theorema}

The results mentioned above make it natural to ask an open question: \emph{if $1<p<2$, whether the results related to Theorems \ref{thm-a} and \ref{thm-b} hold?} An affirmative answer will be given in Theorems \ref{thm-1.1}, \ref{thm-1.2} and \ref{thm-1.3} below.

\subsection{Main results and some related comments}
\label{sect-1.2}

\noindent Before stating our main results, we present the following identities that are essential in analyzing the remainder terms of \eqref{1.6}.

\begin{lem}
[{\rm{\cite[Lemmas 3.3 and 3.4]{Cazacu24}}}]
\label{lem-2.1}
Assume that $1<p<N$. For each magnetic potential $\mathbf{A}$ and complex functions $u\in C_0^\infty(\mathbb{R}^N\setminus\{\mathbf{0}\})$, there holds
\begin{align}\label{2.11}
\mathcal{H}_{\mathbf{A},p}(u)
=\int_{\mathbb{R}^N}
\mathcal{K}_p\left(\nabla\!_\mathbf{A} u, |x|^{\frac{p-N}{p}}\nabla\!_\mathbf{A} \left(u|x|^{\frac{N-p}{p}}\right)\right)\mathrm{d}x,
\end{align}
where
\begin{equation}\label{2.9}
\mathcal{K}_p(\eta,\zeta)
:=\left|\eta\right|^p
-\left|\eta-\zeta\right|^p
-p\left|\eta-\zeta\right|^{p-2}
\mathrm{Re}\left(\eta-\zeta\right)\cdot\bar{\zeta},
\ \
\eta,\zeta\in\mathbb{C}^N.
\end{equation}
Especially, if $\mathbf{A}=0$,
\begin{equation}\label{2.10}
\mathcal{H}_{0,p}(u)
=\int_{\mathbb{R}^N}
\mathcal{K}_p\left(\nabla u,|x|^{\frac{p-N}{p}}\nabla \left(u|x|^{\frac{N-p}{p}}\right)\right)\mathrm{d}x.
\end{equation}
\end{lem}

Studying remainder terms of \eqref{1.6} suffices to consider the right-hand side term of \eqref{2.11}.

\subsubsection{Remainder estimates involving magnetic gradients for the magnetic $L^p$-Hardy inequalities}\label{sect-1.2.1}

\noindent From \cite[Lemma 3.1]{Barbatis04} or \cite[Lemma 4.2]{Lindqvist90}, we know that for $1<p<2$ and $\eta,\zeta\in\mathbb{R}^N$,
\begin{equation*}
\mathcal{K}_p(\eta,\zeta)
\ge c(p)\left(|\eta|+|\eta-\zeta|\right)^{p-2}|\zeta|^2,
\end{equation*}
for some $c(p)>0$. Therefore, inspired by this estimate, we first analyze a family of remainder estimates involving magnetic gradients of \eqref{1.6}.

\begin{thm}\label{thm-1.1}
Let $1<p<2\le N$. Suppose that $\mathbf{B}$ is nontrivial, smooth and closed, then for all $\mathbf{A}$ satisfying $\mathrm{d}\mathbf{A}=\mathbf{B}$,
\begin{align*}
\mathcal{H}_{\mathbf{A},p}(u)
\ge c_1(p)\int_{\mathbb{R}^N}
\left(\left|\nabla\!_\mathbf{A} u\right|
+\frac{N-p}{p}\frac{|u|}{|x|}
\right)^{p-2}
|x|^{\frac{2(p-N)}{p}}
\left|\nabla\!_\mathbf{A} \left(u|x|^{\frac{N-p}{p}}\right)\right|^2
\mathrm{d}x,
\ \
\forall
u\in\mathcal{D}_{\mathbf{A},p},
\end{align*}
where $c_1(p)$ is an explicit constant defined by
\begin{align*}
c_1(p):=\inf_{s^2+t^2>0}
\frac{\left(t^2+s^2+2s+1\right)^{\frac{p}{2}}-1-ps}
{\left(\sqrt{t^2+s^2+2s+1}+1\right)^{p-2}
(t^2+s^2)}
\in\left(0,\frac{p(p-1)}{2^{p-1}}\right].
\end{align*}
Furthermore, the remainder term is optimal due to the fact that
\begin{align*}
\mathcal{H}_{\mathbf{A},p}(u)
\le c_2(p)\int_{\mathbb{R}^N}
\left(\left|\nabla\!_\mathbf{A} u\right|
+\frac{N-p}{p}\frac{|u|}{|x|}
\right)^{p-2}
|x|^{\frac{2(p-N)}{p}}
\left|\nabla\!_\mathbf{A} \left(u|x|^{\frac{N-p}{p}}\right)\right|^2
\mathrm{d}x,
\ \
\forall
u\in\mathcal{D}_{\mathbf{A},p},
\end{align*}
where $c_2(p)$ is an explicit constant defined by
\begin{align*}
c_2(p):=\sup_{s^2+t^2>0}
\frac{\left(t^2+s^2+2s+1\right)^{\frac{p}{2}}-1-ps}
{\left(\sqrt{t^2+s^2+2s+1}+1\right)^{p-2}
(t^2+s^2)}
\in\left[\frac{p}{2^{p-1}},+\infty\right).
\end{align*}
\end{thm}

If $\mathbf{B}=0$, let us choose $\mathbf{A}=0$, then $\nabla\!_\mathbf{A}=\nabla$, and the following corollary holds.

\begin{cor}\label{cor-1.1}
Let $1<p<2\le N$. Then the following inequality
\begin{align*}
\mathcal{H}_{0,p}(u)
\ge c_1(p)\int_{\mathbb{R}^N}
\left(\left|\nabla u\right|
+\frac{N-p}{p}\frac{|u|}{|x|}
\right)^{p-2}
|x|^{\frac{2(p-N)}{p}}
\!\left|\nabla \left(u|x|^{\frac{N-p}{p}}\right)\right|^2
\mathrm{d}x
\end{align*}
holds for all $u\in W^{1,p}(\mathbb{R}^N\setminus\{\mathbf{0}\})$, with  $c_1(p)$ defined by Theorem \ref{thm-1.1}. Furthermore, the remainder term is optimal due to the fact that
\begin{align*}
\mathcal{H}_{0,p}(u)
\le c_2(p)\int_{\mathbb{R}^N}
\left(\left|\nabla u\right|
+\frac{N-p}{p}\frac{|u|}{|x|}
\right)^{p-2}
|x|^{\frac{2(p-N)}{p}}
\!\left|\nabla \left(u|x|^{\frac{N-p}{p}}\right)\right|^2
\mathrm{d}x,
\end{align*}
for all $u\in W^{1,p}(\mathbb{R}^N\setminus\{\mathbf{0}\})$, with  $c_2(p)$ defined by Theorem \ref{thm-1.1}.
\end{cor}

From \cite[Lemma 2.1]{Figalli22}, for $1<p<2$ and $\eta,\zeta\in\mathbb{R}^N$, there exists a constant $c_0(p)>0$ such that
\begin{equation*}
\mathcal{K}_p(\eta,\zeta)
\ge c_0(p)\min\left\{|\zeta|^p,
|\eta-\zeta|^{p-2}|\zeta|^2\right\}.
\end{equation*}
This inspires us to study another improved version of \eqref{1.6} for $1<p<2$.

\begin{thm}\label{thm-1.2}
Let $1<p<2\le N$. Suppose that $\mathbf{B}$ is nontrivial, smooth and closed, then for all $\mathbf{A}$ satisfying $\mathrm{d}\mathbf{A}=\mathbf{B}$,
\begin{align*}
\mathcal{H}_{\mathbf{A},p}(u)
\ge c_3(p)\!\int_{\mathbb{R}^N}
\!\min\left\{\!|x|^{p-N}\!\left|\nabla\!_\mathbf{A} \!\left(u|x|^{\frac{N-p}{p}}\right)\!\right|^p\!,
\left(\!\frac{N\!-\!p}{p}\!\right)^{p-2}\!|u|^{p-2}
|x|^{4-p-\frac{2N}{p}}
\!\left|\nabla\!_\mathbf{A} \!\left(u|x|^{\frac{N-p}{p}}\right)\!\right|^2\!\right\}
\!\mathrm{d}x,
\end{align*}
for all $u\in\mathcal{D}_{\mathbf{A},p}$, where $c_3(p)$ is an explicit constant stated by
\begin{align*}
c_3(p):
=\min\left\{\inf_{s^2+t^2\ge 1}
\frac{\left(t^2+s^2+2s+1\right)^{\frac{p}{2}}-1-ps}
{\left(t^2+s^2\right)^{\frac{p}{2}}},
\inf_{0<s^2+t^2<1}
\frac{\left(t^2+s^2+2s+1\right)^{\frac{p}{2}}-1-ps}
{t^2+s^2}
\right\}
\end{align*}
satisfying
\[
c_3(p)\in\left(0,\frac{p(p-1)}{2}\right].
\]
\end{thm}

A direct consequence of Theorem \ref{thm-1.2} is the following corollary with $\mathbf{A}=0$.

\begin{cor}\label{cor-1.2}
Let $1<p<2\le N$. Then the following inequality
\begin{align*}
\mathcal{H}_{0,p}(u)
\ge c_3(p)\!\int_{\mathbb{R}^N}
\min\left\{|x|^{p-N}\!\left|\nabla \!\left(u|x|^{\frac{N-p}{p}}\right)\!\right|^p,
\left(\!\frac{N\!-\!p}{p}\!\right)^{p-2}\!|u|^{p-2}
|x|^{4-p-\frac{2N}{p}}
\!\left|\nabla \!\left(u|x|^{\frac{N-p}{p}}\right)\!\right|^2\right\}
\mathrm{d}x
\end{align*}
holds for all $u\in W^{1,p}(\mathbb{R}^N\setminus\{\mathbf{0}\})$, where the constant $c_3(p)$ is defined by Theorem \ref{thm-1.2}.
\end{cor}

\begin{rem}
Let $1<p<2\le N$, for $\eta,\zeta\in\mathbb{C}^N$, there holds
\begin{equation}\label{1.4}
\left(|\eta|+|\eta-\zeta|\right)^{p-2}|\zeta|^2
\le\min\left\{|\zeta|^p,|\eta-\zeta|^{p-2}
|\zeta|^2\right\}
\le3^{2-p}\left(|\eta|+|\eta-\zeta|\right)^{p-2}
|\zeta|^2,
\end{equation}
and Lemma \ref{lem-a} will provide a detailed proof of \eqref{1.4} in Appendix \ref{sect-appendix}. Choosing
\[
\eta
=\nabla\!_\mathbf{A}u,
\ \
\zeta=|x|^{\frac{p-N}{p}}\nabla\!_\mathbf{A} \left(u|x|^{\frac{N-p}{p}}\right)
\]
and
\begin{align}\label{2.6}
\eta-\zeta
&=\nabla\!_\mathbf{A}u
-|x|^{\frac{p-N}{p}}\nabla\!_\mathbf{A} \left(u|x|^{\frac{N-p}{p}}\right)
\nonumber\\&=
\nabla u+i\mathbf{A}u
-|x|^{\frac{p-N}{p}}
\left[\nabla\left(u|x|^{\frac{N-p}{p}}\right)
+i\mathbf{A}u|x|^{\frac{N-p}{p}}\right]
\nonumber\\&=
\nabla u
-|x|^{\frac{p-N}{p}}\nabla \left(u|x|^{\frac{N-p}{p}}\right)
\nonumber\\&=
-\frac{N-p}{p}\frac{u}{|x|}\frac{x}{|x|},
\end{align}
in \eqref{1.4}, we see that
\begin{align}\label{1.5}
&\left(\left|\nabla\!_\mathbf{A} u\right|
+\frac{N-p}{p}\frac{|u|}{|x|}
\right)^{p-2}
|x|^{\frac{2(p-N)}{p}}
\left|\nabla\!_\mathbf{A} \left(u|x|^{\frac{N-p}{p}}\right)\right|^2
\nonumber\\&\quad
\le\min\left\{|x|^{p-N}\left|\nabla\!_\mathbf{A} \left(u|x|^{\frac{N-p}{p}}\right)\right|^p,
\left(\frac{N-p}{p}\right)^{p-2}|u|^{p-2}
|x|^{4-p-\frac{2N}{p}}
\left|\nabla\!_\mathbf{A} \left(u|x|^{\frac{N-p}{p}}\right)\right|^2\right\}
\nonumber\\&\quad
\le 3^{2-p}\left(\left|\nabla\!_\mathbf{A} u\right|
+\frac{N-p}{p}\frac{|u|}{|x|}
\right)^{p-2}
|x|^{\frac{2(p-N)}{p}}
\left|\nabla\!_\mathbf{A} \left(u|x|^{\frac{N-p}{p}}\right)\right|^2,
\end{align}
which indicates that the results of Theorem \ref{thm-1.1} is equivalent to those of Theorem \ref{thm-1.2}.
\end{rem}

\begin{rem}
Some discussions about our main results mentioned above are presented below.
\begin{enumerate}
[itemsep=0pt, topsep=2pt, parsep=0pt]

\item[(1)]
By adding remainder terms that involve magnetic gradients, Theorems \ref{thm-1.1} and \ref{thm-1.2} improve the magnetic $L^p$-Hardy inequalities \eqref{1.6} for $1<p<2$, which reply the question presented in
\cite[Remark 1.1]{Cazacu24}. Thus, this paper can be regarded as the supplementary work of \cite{Cazacu24}.

\item[(2)]
For the magnetic-free case $\mathbf{A}=0$, the remainder estimates of Hardy inequality can be deduced from \cite[Lemma 2.2]{Pinchover08}, while such remainder estimates are sightly different from those in Theorems \ref{thm-1.1} and \ref{thm-1.2}.

\item[(3)]
Corollaries \ref{cor-1.1} and \ref{cor-1.2} improve the classical $L^p$-Hardy inequality \eqref{1.2} for $1<p<2$.

\end{enumerate}
\end{rem}

\subsubsection{Remainder estimates involving logarithmic terms for the magnetic $L^p$-Hardy inequalities}\label{sect-1.2.2}

\noindent Now, we turn our attention to establishing another class of remainder terms involving logarithmic terms of \eqref{1.6}.

\begin{thm}\label{thm-1.3}
Let $1<p<2\le N$. Assume that $\mathbf{B}$ is nontrivial, smooth and closed, then for all $\mathbf{A}$ satisfying $\mathrm{d}\mathbf{A}=\mathbf{B}$, there exists a constant $C_{\mathbf{B},p,N}>0$ such that
\begin{equation*}
\mathcal{H}_{\mathbf{A},p}(u)
\ge\frac{c_3(p)}{C_{\mathbf{B},p,N}}
\int_{\mathbb{R}^N}
\frac{|u|^p}{|x|^p\left(1+|\ln|x||^p\right)}
\mathrm{d}x,
\end{equation*}
for all $u\in\mathcal{D}_{\mathbf{A},p}$ with  $(ux\cdot\nabla\!_\mathbf{A} u)=|x||\nabla\!_\mathbf{A} u||u|$, where $c_3(p)$ is defined by Theorem \ref{thm-1.2}.
\end{thm}

\begin{rem}
We present some comments about Theorem \ref{thm-1.3}  below.
\begin{enumerate}
[itemsep=0pt, topsep=2pt, parsep=0pt]

\item[(1)]
The reason of assuming $(ux\cdot\nabla\!_\mathbf{A} u)=|x||\nabla\!_\mathbf{A} u||u|$ is that, in this case,
\begin{align*}
\left(\left|\nabla\!_\mathbf{A} u\right|
+\frac{N-p}{p}\frac{|u|}{|x|}
\right)^{p-2}
|x|^{\frac{2(p-N)}{p}}
\left|\nabla\!_\mathbf{A} \left(u|x|^{\frac{N-p}{p}}\right)\right|^2
=|x|^{p-N}\left|\nabla\!_\mathbf{A} \left(u|x|^{\frac{N-p}{p}}\right)\right|^p,
\end{align*}
then for $1<p<2$,
\begin{align*}
&\min\left\{|x|^{p-N}\left|\nabla\!_\mathbf{A} \left(u|x|^{\frac{N-p}{p}}\right)\right|^p,
\left(\frac{N-p}{p}\right)^{p-2}|u|^{p-2}
|x|^{4-p-\frac{2N}{p}}
\left|\nabla\!_\mathbf{A} \left(u|x|^{\frac{N-p}{p}}\right)\right|^2\right\}
\nonumber\\&\quad
=|x|^{p-N}\left|\nabla\!_\mathbf{A} \left(u|x|^{\frac{N-p}{p}}\right)\right|^p.
\end{align*}
This assumption is natural since it holds for almost all radial functions.

\item[(2)]
Theorem \ref{thm-1.3} provides an improved version of \eqref{1.6} by adding remainder estimates involving logarithmic terms.

\item[(3)]
Theorem \ref{thm-1.3} fills the results of \cite[Theorem 1.3]{Cazacu24} (see also Theorem \ref{thm-a} above) for $1<p<2$.

\item[(4)]
\cite[Proposition 1.2]{Cazacu24} indicates that, for $1<p<N$, the operator
$-\Delta_p
-\left(\frac{N-p}{p}\right)^p
\frac{|\,\cdot\,|^{p-2}}{|x|^p}$ is critical. However, we infer from Theorem \ref{thm-1.3} that, for $1<p<N$, the operator
$-\Delta_{\mathbf{A},p}
-\left(\frac{N-p}{p}\right)^p
\frac{|\,\cdot\,|^{p-2}}{|x|^p}$ becomes subcritical after adding a magnetic field to the $p$-Laplacian.

\end{enumerate}
\end{rem}

\subsection{Applications: remainder estimates of the cylindrical $L^p$-Hardy inequalities for $1<p<2$}\label{sect-1.4}

\noindent Badiale and Tarantello in
\cite{Badiale02} (or Secchi \emph{et al.} in \cite{Secchi03}) studied the following cylindrical $L^p$-Hardy inequalities: let $1<p<k$, $2\le k\le N$, $x=(x',x'')\in\mathbb{R}^{k}\times\mathbb{R}^{N-k}$ and  $u\in C_0^\infty(\mathbb{R}^N\setminus\{x'=\mathbf{0}\})$,
\begin{equation}\label{1.3}
\int_{\mathbb{R}^N}|\nabla u|^p\mathrm{d}x
\ge\left(\frac{k-p}{p}\right)^p
\int_{\mathbb{R}^N}\frac{|u|^p}{|x'|^p}\mathrm{d}x.
\end{equation}
Moreover, the constant $\left(\frac{k-p}{p}\right)^p$ is optimal.

Our purpose in this subsection is to analyze improved versions of \eqref{1.3} by adding remainder terms on the right side of it. We first present the following cylindrical $L^p$-Hardy identities, which are the key point to investigate remainder terms of \eqref{1.3}.

\begin{theorema}
[{\rm{\cite[Theorem 2]{Kalaman24}}}]
Assume that $1<p<\infty$ and $x=(x',x'')\in\mathbb{R}^{k}\times\mathbb{R}^{N-k}$. Then for each complex functions $u\in C_0^\infty(\mathbb{R}^N\setminus\{x'=\mathbf{0}\})$,
\begin{align*}
\int_{\mathbb{R}^N}
\left|\frac{x'}{|x'|}\cdot\nabla\!_k u\right|^p
\mathrm{d}x
-\left|\frac{k-p}{p}\right|^p
\int_{\mathbb{R}^N}\frac{|u|^p}{|x'|^p}\mathrm{d}x
=\int_{\mathbb{R}^N}
\mathcal{K}_p
\left(\frac{x'}{|x'|}\cdot\nabla\!_{k} u, \frac{x'}{|x'|}\cdot\nabla\!_{k} u
+\frac{k-p}{p}\frac{u}{|x'|}\right)
\mathrm{d}x,
\end{align*}
where $\nabla\!_{k}
=\left(\partial_{x_1},\ldots,\partial_{x_{k}}\right)$ and the function $\mathcal{K}_p(\cdot,\cdot)$ is defined by \eqref{2.9}.
\end{theorema}

Based on these identities, remainder estimates of \eqref{1.3} have been proved in \cite{Kalaman24} for  $p\ge2$. Next, we consider the case $1<p<2$. Similar arguments to those of our main results, we obtain the remainder estimates of \eqref{1.3} for $1<p<2$.

\begin{thm}\label{thm-1.5}
Suppose that $1<p<2\le k\le N$ and $x=(x',x'')\in\mathbb{R}^{k}\times\mathbb{R}^{N-k}$. Then, for each complex functions $u\in C_0^\infty
(\mathbb{R}^N\setminus\{x'=\mathbf{0}\})$,

\begin{enumerate}
[itemsep=0pt, topsep=2pt, parsep=0pt]

\item[(1)]
for constants $c_1(p),c_2(p)>0$ defined in Theorem \ref{thm-1.1},
\begin{align*}
&c_2(p)\int_{\mathbb{R}^N}
\left(\left|\frac{x'}{|x'|}\cdot\nabla\!_{k} u\right|
+\frac{k-p}{p}\frac{|u|}{|x'|}
\right)^{p-2}
\left|\frac{x'}{|x'|}\cdot\nabla\!_{k} u
+\frac{k-p}{p}\frac{u}{|x'|}\right|^2
\mathrm{d}x
\\&\quad
\ge\int_{\mathbb{R}^N}
\left|\frac{x'}{|x'|}\cdot\nabla\!_k u\right|^p
\mathrm{d}x
-\left(\frac{k-p}{p}\right)^p
\int_{\mathbb{R}^N}\frac{|u|^p}{|x'|^p}\mathrm{d}x
\\&\quad
\ge c_1(p)\int_{\mathbb{R}^N}
\left(\left|\frac{x'}{|x'|}\cdot\nabla\!_{k} u\right|
+\frac{k-p}{p}\frac{|u|}{|x'|}
\right)^{p-2}
\left|\frac{x'}{|x'|}\cdot\nabla\!_{k} u
+\frac{k-p}{p}\frac{u}{|x'|}\right|^2
\mathrm{d}x;
\end{align*}

\item[(2)]
for constant $c_3(p)>0$ defined in Theorem \ref{thm-1.2},
{\small
\begin{align*}
&\int_{\mathbb{R}^N}
\left|\frac{x'}{|x'|}\cdot\nabla\!_k u\right|^p
\mathrm{d}x
-\left(\frac{k-p}{p}\right)^p
\int_{\mathbb{R}^N}\frac{|u|^p}{|x'|^p}\mathrm{d}x
\\&\quad
\ge c_3(p)\int_{\mathbb{R}^N}
\min\left\{
\!\left|\frac{x'}{|x'|}\cdot\nabla\!_{k} u
+\frac{k-p}{p}\frac{u}{|x'|}\right|^p,
\left(\frac{k-p}{p}\right)^{p-2}
\frac{|u|^{p-2}}{|x'|^{p-2}}
\!\left|\frac{x'}{|x'|}\cdot\nabla\!_{k} u
+\frac{k-p}{p}\frac{u}{|x'|}\right|^2\right\}
\mathrm{d}x.
\end{align*}}
\end{enumerate}
\end{thm}

There holds the following corollary by the Cauchy-Schwarz inequality and $|\nabla\!_ku|\le|\nabla u|$.

\begin{cor}\label{cor-1.5}
Suppose that $1<p<2\le k\le N$ and $x=(x',x'')\in\mathbb{R}^{k}\times\mathbb{R}^{N-k}$. Then, for each complex functions $u\in C_0^\infty
(\mathbb{R}^N\setminus\{x'=\mathbf{0}\})$,

\begin{enumerate}
[itemsep=0pt, topsep=2pt, parsep=0pt]

\item[(1)]
for constant $c_1(p)>0$ defined in Theorem \ref{thm-1.1},
\begin{align*}
&\int_{\mathbb{R}^N}
\left|\nabla u\right|^p
\mathrm{d}x
-\left(\frac{k-p}{p}\right)^p
\int_{\mathbb{R}^N}\frac{|u|^p}{|x'|^p}\mathrm{d}x
\\&\quad
\ge c_1(p)\int_{\mathbb{R}^N}
\left(\left|\frac{x'}{|x'|}\cdot\nabla\!_{k} u\right|
+\frac{k-p}{p}\frac{|u|}{|x'|}
\right)^{p-2}
\left|\frac{x'}{|x'|}\cdot\nabla\!_{k} u
+\frac{k-p}{p}\frac{u}{|x'|}\right|^2
\mathrm{d}x;
\end{align*}

\item[(2)]
for constant $c_3(p)>0$ defined in Theorem \ref{thm-1.2},
{\small
\begin{align*}
&\int_{\mathbb{R}^N}
\left|\nabla u\right|^p
\mathrm{d}x
-\left(\frac{k-p}{p}\right)^p
\int_{\mathbb{R}^N}\frac{|u|^p}{|x'|^p}\mathrm{d}x
\\&\quad
\ge c_3(p)\int_{\mathbb{R}^N}
\min\left\{
\!\left|\frac{x'}{|x'|}\cdot\nabla\!_{k} u
+\frac{k-p}{p}\frac{u}{|x'|}\right|^p,
\left(\frac{k-p}{p}\right)^{p-2}
\frac{|u|^{p-2}}{|x'|^{p-2}}
\!\left|\frac{x'}{|x'|}\cdot\nabla\!_{k} u
+\frac{k-p}{p}\frac{u}{|x'|}\right|^2\right\}
\mathrm{d}x.
\end{align*}}
\end{enumerate}
\end{cor}

\begin{rem}
We present some comments on the above results.
\begin{enumerate}
[itemsep=0pt, topsep=2pt, parsep=0pt]

\item[(1)]
If $k=N$, the results of Corollary \ref{cor-1.5} are equivalent to those of Corollaries \ref{cor-1.1} and \ref{cor-1.2}.

\item[(2)] Similar arguments to those of our main results can also  be used to obtain remainder terms of some other $L^p$-Hardy-type inequalities, such as $L^p$-Hardy inequalities on stratified Lie groups (with the aid of the identities proven in \cite[Theorem 3]{Kalaman24}) and $L^p$-Hardy inequalities on homogeneous Lie groups (due to the identities proven in \cite[Theorem 4]{Kalaman24}) \emph{etc}.

\end{enumerate}
\end{rem}

\subsection{Structure of this paper}\label{sect-1.3}

\begin{itemize}
[itemsep=0pt, topsep=0pt, parsep=0pt]

\item
In Section \ref{sect-2}, we aim to investigate a family of remainder terms involving magnetic gradients of the magnetic $L^p$-Hardy inequalities with $1<p<2$, and prove Theorems \ref{thm-1.1} and \ref{thm-1.2}, Corollaries \ref{cor-1.1} and \ref{cor-1.2}.

\item
In Section \ref{sect-3}, we establish another family of remainder terms involving logarithmic terms of the magnetic $L^p$-Hardy inequalities for $1<p<2$, and prove Theorem \ref{thm-1.3}.

\end{itemize}

\section{Remainder estimates involving magnetic gradients of the magnetic $L^p$-Hardy inequalities: proof of
Theorem \ref{thm-1.1} and Corollary \ref{cor-1.1}, Theorem \ref{thm-1.2} and Corollary \ref{cor-1.2}}
\label{sect-2}

\noindent In this section, we improve \eqref{1.6} when adding remainder terms involving magnetic gradients on the sign-hand side of it, and prove Theorems \ref{thm-1.1} and \ref{thm-1.2}, Corollaries \ref{cor-1.1} and \ref{cor-1.2}.

We first provide a preliminary calculation that will be used frequently.

\begin{lem}\label{lem-2.6}
Assume that $1<p<2$. For all $s^2+t^2>0$,
\begin{equation*}
\left(t^2+s^2+2s+1\right)^{\frac{p}{2}}-1-ps>0.
\end{equation*}
\end{lem}

\begin{proof}[\textbf{Proof}]
We validate this lemma in two cases: $t=0$ and $t\neq0$.

$\diamondsuit$ If $t=0$, then for $s^2>0$ (\emph{i.e.}, $s\in\mathbb{R}\setminus\{0\}$),
\[
\left(t^2+s^2+2s+1\right)^{\frac{p}{2}}-1-ps
=\left(s^2+2s+1\right)^{\frac{p}{2}}-1-ps
:=g_p(s).
\]
Observe that, for all $s\in\mathbb{R}\setminus\{0\}$,
\begin{align}
g'_p(s)
&=p\left[\left(s^2+2s+1\right)^{\frac{p}{2}-1}(s+1)
-1\right],\label{2.12}
\\
g''_p(s)
&=p(p-1)\left(s^2+2s+1\right)^{\frac{p}{2}-1}\ge0.
\label{2.13}
\end{align}
It is easy to verify that $g''_p(s)>0$ for $s\neq-1$. Based on this, $g'_p(-1)=-p<0$ and $g'_p(0)=0$, it follows that $g'_p(s)<0$ for $s\in(-\infty,0)$ and $g'_p(s)>0$ for $s\in(0,+\infty)$. Then,
\[
g_p(s)
>g_p(0)
=0,
\ \
\mathrm{for}
\
\mathrm{all}
\
s\in\mathbb{R}\setminus\{0\}.
\]

$\diamondsuit$ If $t\neq0$, then for $s\in\mathbb{R}$,
\[
\left(t^2+s^2+2s+1\right)^{\frac{p}{2}}-1-ps
>\left(s^2+2s+1\right)^{\frac{p}{2}}-1-ps
=g_p(s).
\]
From \eqref{2.12} and \eqref{2.13}, we obtain $g_p''(s)\ge0$ and $g_p'(0)=0$, then $g'_p(s)\le 0$ for $s\le 0$ and $g'_p(s)\ge 0$ for $s\ge 0$. This together with $g_p(0)=0$, we get $g_p(s)\ge 0$ for all $s\in\mathbb{R}$, as our desired result.
\end{proof}

\subsection{Proof of Theorem \ref{thm-1.1} and Corollary \ref{cor-1.1}}
\label{sect-2.1}

\noindent Before proving Theorem \ref{thm-1.1} and Corollary \ref{cor-1.1}, we prove some estimates about $\mathcal{K}_p$ in the following two lemmas.

\begin{lem}\label{lem-2.2}
Let $1<p<2\le N$, for $\eta,\zeta\in\mathbb{C}^N$, it holds that
\begin{equation*}
\mathcal{K}_p(\eta,\zeta)
\ge c_1(p)\left(|\eta|+|\eta-\zeta|\right)^{p-2}|\zeta|^2,
\end{equation*}
where
\begin{equation*}
c_1(p):
=\inf_{s^2+t^2>0}
\frac{\left(t^2+s^2+2s+1\right)^{\frac{p}{2}}-1-ps}
{\left(
\sqrt{t^2+s^2+2s+1}+1\right)^{p-2}
(t^2+s^2)}.
\end{equation*}
Moreover, $c_1(p)\in\left(0,\frac{p(p-1)}{2^{p-1}}\right]$.
\end{lem}

\begin{proof}[\textbf{Proof}]
For $\eta,\zeta\in\mathbb{R}^N$, the above inequality can be deduced directly from \cite[Lemma 3.1]{Barbatis04} or \cite[Lemma 4.2]{Lindqvist90}. Here we consider $\eta,\zeta\in\mathbb{C}^N$ with optimal constant.

Inspired by \cite{Cazacu24}, let $\eta-\zeta=a+bi$ and $\zeta=c+di$ for $a,b,c,d\in\mathbb{R}^N$,
\begin{align}\label{2.1}
|\eta|^2&=|a+c|^2+|b+d|^2
\nonumber\\&=|a|^2+|b|^2+2(a\cdot c+b\cdot d)+|c|^2+|d|^2,\nonumber\\
|\eta-\zeta|^2&=|a|^2+|b|^2,\nonumber\\
|\zeta|^2&=|c|^2+|d|^2,
\end{align}
then
\begin{align*}
\mathcal{K}_p(\eta,\zeta)
&=\left||a|^2+|b^2|+2(a\cdot c+b\cdot d)+|c|^2+|d|^2\right|^{\frac{p}{2}}
\nonumber\\&\quad
-\left||a|^2+|b|^2\right|^{\frac{p}{2}}
-p\left||a|^2+|b|^2\right|^{\frac{p}{2}-1}
\left(a\cdot c+b\cdot d\right),
\end{align*}
and
\begin{align*}
&\left(|\eta|+|\eta-\zeta|\right)^{p-2}|\zeta|^2
\\&\quad
=\left(
\sqrt{|a|^2+|b|^2+2(a\cdot c+b\cdot d)+|c|^2+|d|^2}
+\sqrt{|a|^2+|b|^2}\right)^{p-2}
\left(|c|^2+|d|^2\right).
\end{align*}

If $|a|^2+|b|^2=0$ or $|c|^2+|d|^2=0$, it is easy to get the desired estimate.

If $|a|^2+|b|^2\neq 0$ and $|c|^2+|d|^2\neq 0$, then it is enough to take into account that
\begin{align*}
c_1(p)
=\inf\limits_{|a|^2+|b|^2\neq 0\atop|c|^2+|d|^2\neq 0}
\frac{\left[
\begin{matrix}
\left||a|^2+|b^2|+2(a\cdot c+b\cdot d)+|c|^2+|d|^2\right|^{\frac{p}{2}}\\[2mm]
-\left||a|^2+|b|^2\right|^{\frac{p}{2}}
-p\left||a|^2+|b|^2\right|^{\frac{p}{2}-1}
\left(a\cdot c+b\cdot d\right)
\end{matrix}
\right]}
{\left(
\sqrt{|a|^2+|b|^2+2(a\cdot c+b\cdot d)+|c|^2+|d|^2}
+\sqrt{|a|^2+|b|^2}\right)^{p-2}
\left(|c|^2+|d|^2\right)}.
\end{align*}
Let $s=\frac{a\cdot c+b\cdot d}{|a|^2+|b|^2}$ and $s^2+t^2=\frac{|c|^2+|d|^2}{|a|^2+|b|^2}$ given that $s^2=\frac{(a\cdot c+b\cdot d)^2}{(|a|^2+|b|^2)^2}
\le\frac{|c|^2+|d|^2}{|a|^2+|b|^2}$. Then, the equality above is reduced to
\begin{equation*}
c_1(p)
=\inf_{s^2+t^2>0}
\frac{\left(t^2+s^2+2s+1\right)^{\frac{p}{2}}-1-ps}
{\left(\sqrt{t^2+s^2+2s+1}+1\right)^{p-2}
(t^2+s^2)}.
\end{equation*}
For simplicity, let
\begin{equation}\label{2.3}
\mathcal{G}(s,t):
=\frac{\left(t^2+s^2+2s+1\right)^{\frac{p}{2}}-1-ps}
{\left(\sqrt{t^2+s^2+2s+1}+1\right)^{p-2}
(t^2+s^2)}.
\end{equation}
To complete the proof of this lemma, it remains to check that
\begin{equation*}
c_1(p)
=\inf_{s^2+t^2>0}\mathcal{G}(s,t)
\in\left(0,\frac{p(p-1)}{2^{p-1}}\right].
\end{equation*}
Actually, choosing $t=0$ and $s\to 0^+$, we arrive at
\begin{align*}
\lim\limits_{t=0\atop s\to 0^+}\mathcal{G}(s,t)
&=\lim\limits_{s\to 0^+}
\frac{\left(s^2+2s+1\right)^{\frac{p}{2}}-1-ps}
{\left(\sqrt{s^2+2s+1}+1\right)^{p-2}
s^2}
\\&=\lim\limits_{s\to 0^+}
\frac{\left(s+1\right)^{p}-1-ps}
{\left(s+2\right)^{p-2}s^2}
\\&=\lim\limits_{s\to 0^+}
\frac{p(s+1)^{p-1}-p}
{(p-2)(s+2)^{p-3}s^2+2s(s+2)^{p-2}}
\\&=\lim\limits_{s\to 0^+}
\frac{p(s+1)^{p-1}-p}
{(ps^2+4s)(s+2)^{p-3}}
\\&=\lim\limits_{s\to 0^+}
\frac{p(p-1)(s+1)^{p-2}}
{2(ps+2)(s+2)^{p-3}+(p-3)(ps^2+4s)(s+2)^{p-4}}
\\&=\frac{p(p-1)}{2^{p-1}},
\end{align*}
which indicates that $c_1(p)\le\frac{p(p-1)}{2^{p-1}}$.

Now, it suffices to verify $c_1(p)>0$. According to Lemma \ref{lem-2.6}, it can be shown that $\mathcal{G}(s,t)>0$ for all $s^2+t^2>0$. As $s^2+t^2\to+\infty$, we get  $\mathcal{G}(s,t)\to1$. As $s^2+t^2\to 0^+$,
\begin{align}\label{2.5}
&\limsup\limits_{s^2+t^2\to 0^+}\mathcal{G}(s,t)
\nonumber\\&\quad=
\limsup\limits_{s^2+t^2\to 0^+}
\frac{\left(t^2+s^2+2s+1\right)^{\frac{p}{2}}-1-ps}
{2^{p-2}(t^2+s^2)}
\nonumber\\&\quad\ge
\inf_{\theta\in [0,\pi]}\lim\limits_{r\to 0^+}
\frac{\left(r^2+2r\cos\theta+1\right)^{\frac{p}{2}}
-1-pr\cos\theta}{2^{p-2}r^2}
\nonumber\\&\quad=
\inf_{\theta\in [0,\pi]}\lim\limits_{r\to 0^+}
\frac{p
\left(r^2+2r\cos\theta+1\right)^{\frac{p}{2}-1}
(r+\cos\theta)-p\cos\theta}{2^{p-1}r}
\nonumber\\&\quad=
\inf_{\theta\in [0,\pi]}\lim\limits_{r\to 0^+}
\frac{p(p-2)
\left(r^2+2r\cos\theta+1\right)^{\frac{p}{2}-2}
(r+\cos\theta)^2
+p\left(r^2+2r\cos\theta+1\right)^{\frac{p}{2}-1}}
{2^{p-1}}
\nonumber\\&\quad=
\inf_{\theta\in[0,\pi]}
\frac{p(p-2)\cos^2\theta+p}{2^{p-1}}
\nonumber\\&\quad=
\frac{p(p-1)}{2^{p-1}}>0,
\end{align}
where $t=r\sin\theta$ and $s=r\cos\theta$. Considering this, $c_1(p)>0$. This completes the proof.
\end{proof}

\begin{lem}\label{lem-2.3}
Let $1<p<2\le N$, for $\eta,\zeta\in\mathbb{C}^N$, there holds
\begin{equation*}
\mathcal{K}_p(\eta,\zeta)\le c_2(p)\left(|\eta|+|\eta-\zeta|\right)^{p-2}|\zeta|^2,
\end{equation*}
where
\begin{equation*}
c_2(p):
=\sup_{s^2+t^2>0}
\frac{\left(t^2+s^2+2s+1\right)^{\frac{p}{2}}-1-ps}
{\left(
\sqrt{t^2+s^2+2s+1}+1\right)^{p-2}
(t^2+s^2)}.
\end{equation*}
Moreover, $c_2(p)\in\left[\frac{p}{2^{p-1}},{+\infty}\right)$.
\end{lem}

\begin{proof}[\textbf{Proof}]
As in the proof of Lemma \ref{lem-2.2}, let $\eta-\zeta=a+bi$ and $\zeta=c+di$ for $a,b,c,d\in\mathbb{R}^N$.
If $|a|^2+|b|^2=0$ or $|c|^2+|d|^2=0$, it is evident that the result is valid. So we will assume that $|a|^2+|b|^2\neq 0$ and $|c|^2+|d|^2\neq 0$, then it suffices to consider that
\begin{align*}
c_2(p)
=\sup\limits_{|a|^2+|b|^2\neq 0\atop|c|^2+|d|^2\neq 0}
\frac{\left[
\begin{matrix}
\left||a|^2+|b^2|+2(a\cdot c+b\cdot d)+|c|^2+|d|^2\right|^{\frac{p}{2}}\\[2mm]
-\left||a|^2+|b|^2\right|^{\frac{p}{2}}
-p\left||a|^2+|b|^2\right|^{\frac{p}{2}-1}
\left(a\cdot c+b\cdot d\right)
\end{matrix}
\right]}
{\left(
\sqrt{|a|^2+|b|^2+2(a\cdot c+b\cdot d)+|c|^2+|d|^2}
+\sqrt{|a|^2+|b|^2}\right)^{p-2}
\left(|c|^2+|d|^2\right)}.
\end{align*}
Let $s=\frac{a\cdot c+b\cdot d}{|a|^2+|b|^2}$ and $s^2+t^2=\frac{|c|^2+|d|^2}{|a|^2+|b|^2}$. Then the above equality turns into
\begin{equation*}
c_2(p)
=\sup_{s^2+t^2>0}
\frac{\left(t^2+s^2+2s+1\right)^{\frac{p}{2}}-1-ps}
{\left(
\sqrt{t^2+s^2+2s+1}+1\right)^{p-2}
(t^2+s^2)}.
\end{equation*}
Thus, it remains to demonstrate that
\begin{equation*}
c_2(p)
=\sup_{s^2+t^2>0}\mathcal{G}(s,t)
\in\left[\frac{p}{2^{p-1}},{+\infty}\right),
\end{equation*}
where $\mathcal{G}(s,t)$ is defined by \eqref{2.3}. Indeed, choosing $s=0$ and $t\to 0^+$, we get
\begin{align*}
\lim\limits_{s=0\atop t\to 0^+}\mathcal{G}(s,t)
&=
\lim\limits_{t\to 0^+}
\frac{\left(t^2+1\right)^{\frac{p}{2}}-1}
{\left(\sqrt{t^2+1}+1\right)^{p-2}t^2}
\\&=
\lim\limits_{t\to 0^+}
\frac{\left(t^2+1\right)^{\frac{p}{2}}-1}
{2^{p-2}t^2}
\\&=
\lim\limits_{t\to 0^+}
\frac{p\left(t^2+1\right)^{\frac{p}{2}-1}t}{2^{p-1}t}
\\&=\frac{p}{2^{p-1}},
\end{align*}
which suggests that $c_2(p)\ge\frac{p}{2^{p-1}}$.

As $s^2+t^2\to+\infty$, we obtain $\mathcal{G}(s,t)\to1$. As $s^2+t^2\to 0^+$,
\begin{align*}
&\liminf\limits_{s^2+t^2\to 0^+}\mathcal{G}(s,t)
\nonumber\\&\quad=
\liminf\limits_{s^2+t^2\to 0^+}
\frac{\left(t^2+s^2+2s+1\right)^{\frac{p}{2}}-1-ps}
{2^{p-2}(t^2+s^2)}
\nonumber\\&\quad\le
\sup_{\theta\in [0,\pi]}\lim\limits_{r\to 0^+}
\frac{\left(r^2+2r\cos\theta+1\right)^{\frac{p}{2}}
-1-pr\cos\theta}{2^{p-2}r^2}
\nonumber\\&\quad=
\sup_{\theta\in[0,\pi]}\lim\limits_{r\to 0^+}
\frac{p
\left(r^2+2r\cos\theta+1\right)^{\frac{p}{2}-1}
(r+\cos\theta)-p\cos\theta}{2^{p-1}r}
\nonumber\\&\quad=
\sup_{\theta\in[0,\pi]}\lim\limits_{r\to 0^+}
\frac{p(p-2)
\left(r^2+2r\cos\theta+1\right)^{\frac{p}{2}-2}
(r+\cos\theta)^2
+p\left(r^2+2r\cos\theta+1\right)^{\frac{p}{2}-1}}
{2^{p-1}}
\nonumber\\&\quad=
\sup_{\theta\in[0,\pi]}
\frac{p(p-2)\cos^2\theta+p}{2^{p-1}}
\nonumber\\&\quad=
\frac{p}{2^{p-1}},
\end{align*}
this, along with \eqref{2.5}, indicates that $0<\lim_{s^2+t^2\to 0}\mathcal{G}(s,t)<+\infty$. The proof is completed.
\end{proof}

Now, we are in a position to prove Theorem \ref{thm-1.1} and Corollary \ref{cor-1.1}.

\begin{proof}[\textbf{Proof of Theorem \ref{thm-1.1}}]
One hand, from \eqref{2.11}, Lemma \ref{lem-2.2} and \eqref{2.6}, we derive that
\begin{align*}
\mathcal{H}_{\mathbf{A},p}(u)
&\ge c_1(p)\int_{\mathbb{R}^N}
\left[\left|\nabla\!_\mathbf{A} u\right|
+\left|\nabla\!_\mathbf{A} u
-|x|^{\frac{p-N}{p}}\nabla\!_\mathbf{A} \left(u|x|^{\frac{N-p}{p}}\right)\right|
\right]^{p-2}
|x|^{\frac{2(p-N)}{p}}
\left|\nabla\!_\mathbf{A} \left(u|x|^{\frac{N-p}{p}}\right)\right|^2
\mathrm{d}x
\\&
=c_1(p)\int_{\mathbb{R}^N}
\left(\left|\nabla\!_\mathbf{A} u\right|
+\frac{N-p}{p}\frac{|u|}{|x|}
\right)^{p-2}
|x|^{\frac{2(p-N)}{p}}
\left|\nabla\!_\mathbf{A} \left(u|x|^{\frac{N-p}{p}}\right)\right|^2
\mathrm{d}x,
\end{align*}
where $c_1(p)$ is a constant defined by Lemma \ref{lem-2.2} satisfying $c_1(p)\in\left(0,\frac{p(p-1)}{2^{p-1}}\right]$.

On the other hand, it can be seen from \eqref{2.11}, Lemma \ref{lem-2.3} and \eqref{2.6} that
\begin{align*}
\mathcal{H}_{\mathbf{A},p}(u)
&\le
c_2(p)\int_{\mathbb{R}^N}
\left[\left|\nabla\!_\mathbf{A} u\right|
+\left|\nabla\!_\mathbf{A} u
-|x|^{\frac{p-N}{p}}\nabla\!_\mathbf{A} \left(u|x|^{\frac{N-p}{p}}\right)\right|
\right]^{p-2}
|x|^{\frac{2(p-N)}{p}}
\left|\nabla\!_\mathbf{A} \left(u|x|^{\frac{N-p}{p}}\right)\right|^2
\mathrm{d}x
\\&=
c_2(p)\int_{\mathbb{R}^N}
\left(\left|\nabla\!_\mathbf{A} u\right|
+\frac{N-p}{p}\frac{|u|}{|x|}
\right)^{p-2}
|x|^{\frac{2(p-N)}{p}}
\left|\nabla\!_\mathbf{A} \left(u|x|^{\frac{N-p}{p}}\right)\right|^2
\mathrm{d}x,
\end{align*}
where $c_2(p)\in\left[\frac{p}{2^{p-1}},{+\infty}\right)$ is a constant defined by Lemma \ref{lem-2.3}. This completes the proof.
\end{proof}

\begin{proof}
[\textbf{Proof of Corollary \ref{cor-1.1}}]
By substituting \eqref{2.11} for \eqref{2.10} and using similar arguments to those of the proof of Theorem \ref{thm-1.1}, we can complete the proof of Corollary \ref{cor-1.1}.
\end{proof}

\subsection{Proof of Theorem \ref{thm-1.2} and Corollary \ref{cor-1.2}}
\label{sect-2.2}

\noindent In this subsection, we first establish a crucial lemma (\emph{i.e.}, Lemma \ref{lem-2.4}) that enables us to prove Theorem \ref{thm-1.2} and Corollary \ref{cor-1.2}.

\begin{lem}\label{lem-2.4}
Let $1<p<2\le N$, for $\eta,\zeta\in\mathbb{C}^N$, we get
\begin{equation*}
\mathcal{K}_p(\eta,\zeta)
\ge c_3(p)\min\left\{|\zeta|^p,
|\eta-\zeta|^{p-2}|\zeta|^2\right\},
\end{equation*}
where
\begin{equation*}
c_3(p):
=\min\left\{\inf_{s^2+t^2\ge 1}
\frac{\left(t^2+s^2+2s+1\right)
^{\frac{p}{2}}-1-ps}
{\left(t^2+s^2\right)^{\frac{p}{2}}},
\inf_{0<s^2+t^2<1}
\frac{\left(t^2+s^2+2s+1\right)
^{\frac{p}{2}}-1-ps}
{t^2+s^2}
\right\}.
\end{equation*}
Moreover, $c_3(p)\in\left(0,\frac{p(p-1)}{2}\right]$.
\end{lem}

\begin{proof}[\textbf{Proof}]
For $\eta,\zeta\in\mathbb{R}^N$, the above inequality can be deduced directly from \cite[Lemma 2.1]{Figalli22} but with no explicit constant. Here we consider $\eta,\zeta\in\mathbb{C}^N$ with optimal constant.

Following Lemma \ref{lem-2.2}, let $\eta-\zeta=a+bi$ and $\zeta=c+di$ for $a,b,c,d\in\mathbb{R}^N$,
 and the transformations as \eqref{2.1} give
\begin{align*}
\min\left\{|\zeta|^p,|\eta-\zeta|^{p-2}|\zeta|^2\right\}
=\min\left\{\big||c|^2+|d|^2\big|^{\frac{p}{2}},
\big||a|^2+|b|^2\big|^{\frac{p}{2}-1}
\big(|c|^2+|d|^2\big)\right\}.
\end{align*}
It is obvious that the result holds when $|a|^2+|b|^2=0$ or $|c|^2+|d|^2=0$. So we will assume that $|a|^2+|b|^2\neq 0$ and $|c|^2+|d|^2\neq 0$, it is sufficient to consider that
\begin{equation*}
c_3(p)
=\inf\limits_{|a|^2+|b|^2\neq 0\atop|c|^2+|d|^2\neq 0}
\frac{\left[
\begin{matrix}
\left||a|^2+|b^2|+2(a\cdot c+b\cdot d)+|c|^2+|d|^2\right|^{\frac{p}{2}}\\[2mm]
-\left||a|^2+|b|^2\right|^{\frac{p}{2}}
-p\left||a|^2+|b|^2\right|^{\frac{p}{2}-1}(a\cdot c+b\cdot d)
\end{matrix}
\right]}
{\min\left\{\big||c|^2+|d|^2\big|^{\frac{p}{2}},
\big||a|^2+|b|^2\big|^{\frac{p}{2}-1}
\big(|c|^2+|d|^2\big)\right\}}.
\end{equation*}
Let $s=\frac{a\cdot c+b\cdot d}{|a|^2+|b|^2}$ and $s^2+t^2=\frac{|c|^2+|d|^2}{|a|^2+|b|^2}$. Therefore, the above equality changes into
\begin{equation*}
c_3(p)
=\inf_{s^2+t^2>0}
\frac{\left(t^2+s^2+2s+1\right)^{\frac{p}{2}}-1-ps}
{\min\left\{\left(t^2+s^2\right)^{\frac{p}{2}},
t^2+s^2\right\}}
\in\left(0,\frac{p(p-1)}{2}\right].
\end{equation*}
Indeed, choosing $t=0$ and $s\to 0^+$, we deduce that \begin{align*}
\lim\limits_{t=0\atop s\to 0^+}
\frac{\left(t^2+s^2+2s+1\right)^{\frac{p}{2}}-1-ps}
{\min\left\{\left(t^2+s^2\right)^{\frac{p}{2}},
t^2+s^2\right\}}
&=
\lim\limits_{s\to 0^+}
\frac{\left(s^2+2s+1\right)^{\frac{p}{2}}-1-ps}{s^2}
\\&=
\lim\limits_{s\to 0^+}
\frac{\left(s+1\right)^p-1-ps}{s^2}
\\&=
\lim\limits_{s\to 0^+}
\frac{p\left(s+1\right)^{p-1}-p}{2s}
\\&=
\lim\limits_{s\to 0^+}
\frac{p(p-1)\left(s+1\right)^{p-2}}{2}
\\&=
\frac{p(p-1)}{2},
\end{align*}
this means that $c_3(p)\le\frac{p(p-1)}{2}$.

In the following, we will show that $c_3(p)>0$. We proceed with the proof in two cases.

-- \textbf{Case $1$: $s^2+t^2\ge 1$.} In this case, it is obvious that
\[
\frac{\left(t^2+s^2+2s+1\right)^{\frac{p}{2}}-1-ps}
{\min\left\{\left(t^2+s^2\right)^{\frac{p}{2}},
t^2+s^2\right\}}
=\frac{\left(t^2+s^2+2s+1\right)^{\frac{p}{2}}-1-ps}
{\left(t^2+s^2\right)^{\frac{p}{2}}}
:=\mathcal{F}_1(s,t).
\]
Lemma \ref{lem-2.6} allows us to prove that $\mathcal{F}_1(s,t)>0$ for all $s^2+t^2\ge 1$.

As $s^2+t^2\to+\infty$, we obtain $\mathcal{F}_1(s,t)\to 1$. When $s^2+t^2=1$, then $s\in [-1,1]$ and $\mathcal{F}_1(s,t)
=2^{\frac{p}{2}}\left(s+1\right)^{\frac{p}{2}}-1-ps:
=h_p(s)$. Notice that
$h'_p(s)
=p\left[2^{\frac{p}{2}-1}
\left(s+1\right)^{\frac{p}{2}-1}-1\right]$. Let $h'_p(s_*)=0$, then $s_*=-\frac{1}{2}$, and $h_p(s_*)=\frac{p}{2}>0$ thanks to $p\in(1,2)$. Furthermore, $h'_p(s)\ge 0$ for $s\in[-1,s_*]$, and $h'_p(s)\le0$ for $s\in[s_*,1]$. Hence, we deduce that
\[
\mathcal{F}_1(s,t)
=h_p(s)
\ge\min\left\{h_p(-1),h_p(1)\right\}
=\min\left\{p-1,2^p-p-1\right\}
>0,
\]
when $s^2+t^2=1$.

-- \textbf{Case $2$: $0<s^2+t^2<1$.} In this case, we get
\[
\frac{\left(t^2+s^2+2s+1\right)^{\frac{p}{2}}-1-ps}
{\min\left\{\left(t^2+s^2\right)^{\frac{p}{2}},
t^2+s^2\right\}}
=\frac{\left(t^2+s^2+2s+1\right)^{\frac{p}{2}}-1-ps}
{t^2+s^2}
:=\mathcal{F}_2(s,t).
\]
It deduces from Lemma \ref{lem-2.6} that $\mathcal{F}_2(s,t)>0$ for all $0<s^2+t^2<1$.

When $s^2+t^2\to 1^-$, $\mathcal{F}_2(s,t)
\sim 2^{\frac{p}{2}}\left(s+1\right)^{\frac{p}{2}}-1-ps$, then similar to the Case $1$ when $s^2+t^2=1$, we have $\lim_{s^2+t^2\to1^-}\mathcal{F}_2(s,t)>0$. Now, we consider $s^2+t^2\to 0^+$. Observe that
\begin{align*}
\liminf_{s^2+t^2\to 0^+}
\frac{\left(t^2+s^2+2s+1\right)^{\frac{p}{2}}-1-ps}
{\frac{p(p-1)}{2}(t^2+s^2)}
&\ge
\liminf_{s^2+t^2\to 0^+}
\frac{\left(t^2+s^2+2s+1\right)^{\frac{p}{2}}}
{\frac{p(p-1)}{2}(t^2+s^2)+ps+1}
\\&\ge
\liminf_{s^2+t^2\to 0^+}
\frac{\left(t^2+s^2+2s+1\right)^{\frac{p}{2}}}
{\frac{p}{2}(t^2+s^2)+ps+1}
\\&=
\liminf_{s^2+t^2\to 0^+}
\frac{\left(t^2+s^2+2s+1\right)^{\frac{p}{2}}}
{\frac{p}{2}(t^2+s^2+2s+1)+\frac{2-p}{2}}
\\&\ge
\inf_{\theta\in[0,\pi]}\lim_{r\to 0^+}
\frac{\left(r^2+2r\cos\theta+1\right)^{\frac{p}{2}}}
{\frac{p}{2}(r^2+2r\cos\theta+1)+\frac{2-p}{2}}
\\&=
1,
\end{align*}
where $t=r\sin\theta$ and $s=r\cos\theta$, and the first inequality follows from $\lim_{n\to\infty}\frac{a_n}{b_n}\ge \lim_{n\to\infty}\frac{a_n+c_n}{b_n+c_n}$, if $\lim_{n\to\infty}\frac{a_n+c_n}{b_n+c_n}\ge 1$, $b_n>0$ and $c_n>0$. Therefore,
\begin{equation*}
\liminf_{s^2+t^2\to 0^+}\mathcal{F}_2(s,t)
=\liminf_{s^2+t^2\to 0^+}
\frac{\left(t^2+s^2+2s+1\right)^{\frac{p}{2}}-1-ps}
{t^2+s^2}
\ge\frac{p(p-1)}{2}
>0.
\end{equation*}
To summarize, Case $1$ and Case $2$ reveal that $c_3(p)>0$. This ends the proof of Lemma \ref{lem-2.4}.
\end{proof}

Now, we are ready to prove Theorem \ref{thm-1.2} and Corollary \ref{cor-1.2}.

\begin{proof}[\textbf{Proof of Theorem \ref{thm-1.2}}]
From \eqref{2.11}, Lemma \ref{lem-2.4} and \eqref{2.6}, we get
{\small
\begin{align*}
&\mathcal{H}_{\mathbf{A},p}(u)
\\&\ \
\ge c_3(p)\int_{\mathbb{R}^N}
\min\left\{|x|^{p-N}
\!\left|
\nabla\!_\mathbf{A}\!\left(u|x|^{\frac{N-p}{p}}\right)\!
\right|^p,
\left|
\nabla\!_\mathbf{A}u
\!-\!|x|^{\frac{p-N}{p}}\nabla\!_\mathbf{A}
\!\left(u|x|^{\frac{N-p}{p}}\right)\!
\right|^{p-2}
|x|^{\frac{2(p-N)}{p}}
\!\left|
\nabla\!_\mathbf{A}\!\left(u|x|^{\frac{N-p}{p}}\right)\!
\right|^2
\right\}\mathrm{d}x
\\&\ \
=c_3(p)\int_{\mathbb{R}^N}
\min\left\{|x|^{p-N}\left|\nabla\!_\mathbf{A} \left(u|x|^{\frac{N-p}{p}}\right)\right|^p,
\left(\frac{N-p}{p}\right)^{p-2}|u|^{p-2}
|x|^{4-p-\frac{2N}{p}}
\left|\nabla\!_\mathbf{A} \left(u|x|^{\frac{N-p}{p}}\right)\right|^2\right\}
\mathrm{d}x.
\end{align*}}\noindent
The proof of Theorem \ref{thm-1.2} is completed.
\end{proof}

\begin{proof}
[\textbf{Proof of Corollary \ref{cor-1.2}}]
By replacing \eqref{2.11} with \eqref{2.10} and using arguments analogous to those used in the proof of Theorem \ref{thm-1.2}, we finish the proof of Corollary \ref{cor-1.2}.
\end{proof}

\section{Remainder estimates involving logarithmic terms of the magnetic $L^p$-Hardy inequalities: proof of Theorem \ref{thm-1.3}}\label{sect-3}

\noindent In this section, we focus on improving \eqref{1.6} by adding remainder estimates involving logarithmic terms, and prove Theorem \ref{thm-1.3}. We first present some preliminary lemmas derived in \cite{Cazacu24}.

\begin{lem}
[{\rm{\cite[Lemma 3.1]{Cazacu24}}}]
\label{lem-3.1}
Let $1<p<N$. For all $\widehat{R}>0$, let us denote $B_{\widehat{R}}(0):
=\big\{x\in\mathbb{R}^N:
|x|<\widehat{R}\big\}$ and
$B_{\widehat{R}}^c(0):
=\big\{x\in\mathbb{R}^N:
|x|\ge\widehat{R}\big\}$, then
\begin{enumerate}
[itemsep=0pt, topsep=2pt, parsep=0pt]

\item[(1)]
for all $u\in C_0^\infty(B_{\widehat{R}}(0))$,
\begin{align*}
\int_{B_{\widehat{R}}(0)}
\frac{\left|u\right|^p}
{|x|^p\left|\ln\frac{{\widehat{R}}}{|x|}\right|^{p}}
\mathrm{d}x
\le \left(\frac{p}{p-1}\right)^p
\int_{B_{\widehat{R}}(0)}|x|^{p-N}\left|\nabla \left(u|x|^{\frac{N-p}{p}}\right)\right|^p\mathrm{d}x;
\end{align*}

\item[(2)]
for all $u\in C_0^\infty(B_{\widehat{R}}^c(0))$,
\begin{align*}
\int_{B_{\widehat{R}}^c(0)}
\frac{\left|u\right|^p}
{|x|^p
\left|\ln\frac{{\widehat{R}}}{|x|}\right|^{p}}
\mathrm{d}x
\le\left(\frac{p}{p-1}\right)^p
\int_{B_{\widehat{R}}^c(0)}
|x|^{p-N}\left|\nabla \left(u|x|^{\frac{N-p}{p}}\right)\right|^p\mathrm{d}x.
\end{align*}
\end{enumerate}
\end{lem}

\begin{lem}
[{\rm{\cite[Lemma 3.2]{Cazacu24}}}]
\label{lem-3.2}
Let $1<p<\infty$ and $N\ge2$. Assume that $\mathbf{B}\neq0$ and $\mathbf{A}$ is a vector field satisfying $\mathrm{d}\mathbf{A}=\mathbf{B}$. Then there exists a constant $R>1$ such that
\begin{equation*}
\mu_\mathbf{B}(R)
:=
\inf\limits_{u\in W^{1,p}(\mathcal{O}_R)
\atop u\neq0}
\frac{\int_{\mathcal{O}_R}
\left|\nabla\!_\mathbf{A} u\right|^p
\mathrm{d}x}
{\int_{\mathcal{O}_R}|u|^p\mathrm{d}x}
\neq0,
\end{equation*}
where $\mathcal{O}_R:=B_R(0)\setminus B_{\frac{1}{R}}(0)$.
\end{lem}

For $1<p<2$, a more accurate estimate stated as follows is needed.

\begin{lem}\label{lem-3.3}
Let $1<p<2\le N$. Suppose that $\mathbf{B}$ is a nontrivial magnetic field and $\mathbf{A}$ is a vector field satisfying $\mathrm{d}\mathbf{A}=\mathbf{B}$.
Let $R>1$ be obtained in Lemma \ref{lem-3.2}, and define
\begin{equation*}
\nu_\mathbf{B}(R)
:=\inf\limits_{u\in W^{1,p}(\mathcal{O}_R)
\atop u\neq0}
\frac{\int_{\mathcal{O}_R}
\min\left\{\left|\nabla\!_\mathbf{A} u\right|^p,
\left(\frac{N-p}{p}\right)^{p-2}
|u|^{p-2}\left|\nabla\!_\mathbf{A} u\right|^2\right\}
\mathrm{d}x}
{\int_{\mathcal{O}_R}|u|^p\mathrm{d}x},
\end{equation*}
then $\nu_\mathbf{B}(R)\in(0,\infty)$.
\end{lem}

\begin{proof}[\textbf{Proof}]
Choosing
\[
\zeta
=\nabla\!_\mathbf{A}u,
\ \
\eta=|x|^{\frac{p-N}{p}}\nabla\!_\mathbf{A} \left(u|x|^{\frac{N-p}{p}}\right)
\]
and
\begin{equation*}
|\eta-\zeta|
=\left|\nabla\!_\mathbf{A}u
-|x|^{\frac{p-N}{p}}\nabla\!_\mathbf{A} \left(u|x|^{\frac{N-p}{p}}\right)\right|
\xlongequal{\eqref{2.6}}\frac{N-p}{p}\frac{|u|}{|x|}
\end{equation*}
in \eqref{1.4}, we obtain
\begin{align*}
&\min\left\{\left|\nabla\!_\mathbf{A} u\right|^p,
\left(\frac{N-p}{p}\right)^{p-2}
\frac{|u|^{p-2}}{|x|^{p-2}}\left|\nabla\!_\mathbf{A} u\right|^2\right\}
\\&\quad\ge\left[|x|^{\frac{p-N}{p}}
\left|\nabla\!_\mathbf{A} \left(u|x|^{\frac{N-p}{p}}\right)\right|
+\frac{N-p}{p}\frac{|u|}{|x|}\right]^{p-2}
\left|\nabla\!_\mathbf{A} u\right|^2.
\end{align*}
From this, it follows that
\begin{align}\label{3.15}
&\int_{\mathcal{O}_R}
\min\left\{\left|\nabla\!_\mathbf{A} u\right|^p,
\left(\frac{N-p}{p}\right)^{p-2}
|u|^{p-2}\left|\nabla\!_\mathbf{A} u\right|^2\right\}
\mathrm{d}x
\nonumber\\&\quad\ge
R^{p-2}\int_{\mathcal{O}_R}
\min\left\{\left|\nabla\!_\mathbf{A} u\right|^p,
\left(\frac{N-p}{p}\right)^{p-2}
\frac{|u|^{p-2}}{|x|^{p-2}}\left|\nabla\!_\mathbf{A} u\right|^2\right\}
\mathrm{d}x
\nonumber\\&\quad\ge
R^{p-2}\int_{\mathcal{O}_R}
\left[|x|^{\frac{p-N}{p}}
\left|\nabla\!_\mathbf{A} \left(u|x|^{\frac{N-p}{p}}\right)\right|
+\frac{N-p}{p}\frac{|u|}{|x|}\right]^{p-2}
\left|\nabla\!_\mathbf{A} u\right|^2
\mathrm{d}x
\nonumber\\&\quad\ge
R^{p-2}\left\{\int_{\mathcal{O}_R}
\left[|x|^{\frac{p-N}{p}}
\left|\nabla\!_\mathbf{A} \left(u|x|^{\frac{N-p}{p}}\right)\right|
+\frac{N-p}{p}\frac{|u|}{|x|}\right]^p
\mathrm{d}x\right\}^{\frac{p-2}{p}}
\left(\int_{\mathcal{O}_R}\left|\nabla\!_\mathbf{A} u\right|^p
\mathrm{d}x\right)^{\frac{2}{p}}
\nonumber\\&\quad=
R^{p-2}
\left\{
\frac{\int_{\mathcal{O}_R}\left|\nabla\!_\mathbf{A} u\right|^p
\mathrm{d}x}
{\int_{\mathcal{O}_R}
\left[|x|^{\frac{p-N}{p}}
\left|\nabla\!_\mathbf{A} \left(u|x|^{\frac{N-p}{p}}\right)\right|
+\frac{N-p}{p}\frac{|u|}{|x|}\right]^p
\mathrm{d}x}\right\}^{\frac{2}{p}}
\nonumber\\&\qquad\times
\int_{\mathcal{O}_R}
\left[|x|^{\frac{p-N}{p}}
\left|\nabla\!_\mathbf{A} \left(u|x|^{\frac{N-p}{p}}\right)\right|
+\frac{N-p}{p}\frac{|u|}{|x|}\right]^p
\mathrm{d}x
\nonumber\\&\quad\ge
\left(\frac{N-p}{p}\right)^p
R^{-2}
\left\{
\frac{\int_{\mathcal{O}_R}\left|\nabla\!_\mathbf{A} u\right|^p
\mathrm{d}x}
{\int_{\mathcal{O}_R}
\left[|x|^{\frac{p-N}{p}}
\left|\nabla\!_\mathbf{A} \left(u|x|^{\frac{N-p}{p}}\right)\right|
+\frac{N-p}{p}\frac{|u|}{|x|}\right]^p
\mathrm{d}x}\right\}^{\frac{2}{p}}
\int_{\mathcal{O}_R}|u|^p\mathrm{d}x,
\end{align}
where the third inequality derives from the following inequality
{\footnotesize
\begin{align*}
&\left(\int_{\mathcal{O}_R}\left|\nabla\!_\mathbf{A} u\right|^p\mathrm{d}x\right)^{\frac{2}{p}}
\nonumber\\&\quad=
\left\{\int_{\mathcal{O}_R}
\frac{\left|\nabla\!_\mathbf{A} u\right|^p}
{\left[|x|^{\frac{p-N}{p}}
\left|\nabla\!_\mathbf{A} \left(u|x|^{\frac{N-p}{p}}\right)\right|
+\frac{N-p}{p}\frac{|u|}{|x|}\right]^{\frac{p(2-p)}{2}}}
\left[|x|^{\frac{p-N}{p}}
\left|\nabla\!_\mathbf{A} \left(u|x|^{\frac{N-p}{p}}\right)\right|
+\frac{N-p}{p}\frac{|u|}{|x|}\right]^{\frac{p(2-p)}{2}}
\mathrm{d}x\right\}^{\frac{2}{p}}
\nonumber\\&\quad\le
\int_{\mathcal{O}_R}
\frac{\left|\nabla\!_\mathbf{A} u\right|^2}
{\left[|x|^{\frac{p-N}{p}}
\left|\nabla\!_\mathbf{A} \left(u|x|^{\frac{N-p}{p}}\right)\right|
+\frac{N-p}{p}\frac{|u|}{|x|}\right]^{2-p}}
\mathrm{d}x
\left\{\int_{\mathcal{O}_R}
\left[|x|^{\frac{p-N}{p}}
\left|\nabla\!_\mathbf{A} \left(u|x|^{\frac{N-p}{p}}\right)\right|
+\frac{N-p}{p}\frac{|u|}{|x|}\right]^{p}
\mathrm{d}x\right\}^{\frac{2-p}{p}}.
\end{align*}}\noindent
Observe that \eqref{3.15} is equivalent to
\begin{align}\label{3.12}
&\frac{\int_{\mathcal{O}_R}
\min\left\{\left|\nabla\!_\mathbf{A} u\right|^p,
\left(\frac{N-p}{p}\right)^{p-2}
|u|^{p-2}\left|\nabla\!_\mathbf{A} u\right|^2\right\}
\mathrm{d}x}
{\int_{\mathcal{O}_R}|u|^p\mathrm{d}x}
\nonumber\\&\quad\ge
\left(\frac{N-p}{p}\right)^p
R^{-2}
\left\{
\frac{\int_{\mathcal{O}_R}\left|\nabla\!_\mathbf{A} u\right|^p
\mathrm{d}x}
{\int_{\mathcal{O}_R}
\left[|x|^{\frac{p-N}{p}}
\left|\nabla\!_\mathbf{A} \left(u|x|^{\frac{N-p}{p}}\right)\right|
+\frac{N-p}{p}\frac{|u|}{|x|}\right]^p
\mathrm{d}x}\right\}^{\frac{2}{p}}.
\end{align}
With the help of the triangle inequality and  \eqref{2.6}, we get
\begin{align*}
|x|^{\frac{p-N}{p}}
\left|\nabla\!_\mathbf{A} \left(u|x|^{\frac{N-p}{p}}\right)\right|
&=
\left||x|^{\frac{p-N}{p}}
\nabla\!_\mathbf{A} \left(u|x|^{\frac{N-p}{p}}\right)
-\nabla\!_\mathbf{A}u
+\nabla\!_\mathbf{A}u\right|
\\&\le
\left||x|^{\frac{p-N}{p}}
\nabla\!_\mathbf{A} \left(u|x|^{\frac{N-p}{p}}\right)
-\nabla\!_\mathbf{A}u\right|
+\left|\nabla\!_\mathbf{A}u\right|
\\&=%3
\frac{N-p}{p}\frac{|u|}{|x|}
+\left|\nabla\!_\mathbf{A}u\right|,
\end{align*}
from this and Lemma \ref{lem-3.2}, it follows that
\begin{align*}
&\int_{\mathcal{O}_R}
\left[|x|^{\frac{p-N}{p}}
\left|\nabla\!_\mathbf{A} \left(u|x|^{\frac{N-p}{p}}\right)\right|
+\frac{N-p}{p}\frac{|u|}{|x|}\right]^p
\mathrm{d}x
\nonumber\\&\quad\le
\int_{\mathcal{O}_R}
\left[\left|\nabla\!_\mathbf{A} u\right|
+\frac{2(N-p)}{p}\frac{|u|}{|x|}\right]^p
\mathrm{d}x
\nonumber\\&\quad\le
2^{p-1}\left\{\int_{\mathcal{O}_R}
\left|\nabla\!_\mathbf{A} u\right|^p\mathrm{d}x
+\left[\frac{2(N-p)}{p}\right]^p
\int_{\mathcal{O}_R}\frac{|u|^p}{|x|^p}
\mathrm{d}x\right\}
\nonumber\\&\quad\le
2^{p-1}\left\{\int_{\mathcal{O}_R}
\left|\nabla\!_\mathbf{A} u\right|^p\mathrm{d}x
+\left[\frac{2(N-p)}{p}\right]^pR^p
\int_{\mathcal{O}_R}|u|^p
\mathrm{d}x\right\}
\nonumber\\&\quad\le
2^{p-1}
\max\left\{\left[\frac{2(N-p)}{p}\right]^pR^p,1\right\}
\left(\int_{\mathcal{O}_R}
\left|\nabla\!_\mathbf{A} u\right|^p\mathrm{d}x
+\int_{\mathcal{O}_R}|u|^p
\mathrm{d}x\right)
\nonumber\\&\quad\le
2^{p-1}
\max\left\{\left[\frac{2(N-p)}{p}\right]^pR^p,1\right\}
\left[1+\frac{1}{\mu_\mathbf{B}(R)}\right]
\int_{\mathcal{O}_R}
\left|\nabla\!_\mathbf{A} u\right|^p\mathrm{d}x,
\end{align*}
that is,
\begin{align}\label{3.13}
&\frac{\int_{\mathcal{O}_R}
\left|\nabla\!_\mathbf{A} u\right|^p\mathrm{d}x}
{\int_{\mathcal{O}_R}
\left[|x|^{\frac{p-N}{p}}
\left|\nabla\!_\mathbf{A} \left(u|x|^{\frac{N-p}{p}}\right)\right|
+\frac{N-p}{p}\frac{|u|}{|x|}\right]^p
\mathrm{d}x}
\nonumber\\&\quad\ge\frac{1}{2^{p-1}
\max\left\{\left[\frac{2(N-p)}{p}\right]^pR^p,
1\right\}
\left[1+\frac{1}{\mu_\mathbf{B}(R)}\right]}.
\end{align}
Consequently, combining \eqref{3.13} with \eqref{3.12}, it gives that
\begin{align*}
&\frac{\int_{\mathcal{O}_R}
\min\left\{\left|\nabla\!_\mathbf{A} u\right|^p,
\left(\frac{N-p}{p}\right)^{p-2}
|u|^{p-2}\left|\nabla\!_\mathbf{A} u\right|^2\right\}
\mathrm{d}x}
{\int_{\mathcal{O}_R}|u|^p\mathrm{d}x}
\nonumber\\&\quad\ge
\left(\frac{N-p}{p}\right)^p
R^{-2}
\left\{2^{p-1}
\max\left\{\left[\frac{2(N-p)}{p}\right]^pR^p,1\right\}
\left[1+\frac{1}{\mu_\mathbf{B}(R)}\right]
\right\}^{-\frac{2}{p}}>0.
\end{align*}
Furthermore, using
\begin{equation*}
\frac{\int_{\mathcal{O}_R}
\min\left\{\left|\nabla\!_\mathbf{A} u\right|^p,
\left(\frac{N-p}{p}\right)^{p-2}
|u|^{p-2}\left|\nabla\!_\mathbf{A} u\right|^2\right\}
\mathrm{d}x}
{\int_{\mathcal{O}_R}|u|^p\mathrm{d}x}
\le\frac{\int_{\mathcal{O}_R}
\left|\nabla\!_\mathbf{A} u\right|^p
\mathrm{d}x}
{\int_{\mathcal{O}_R}|u|^p\mathrm{d}x},
\end{equation*}
we conclude that
\[
\left(\frac{N-p}{p}\right)^p
R^{-2}
\left\{2^{p-1}
\max\left\{\left[\frac{2(N-p)}{p}\right]^pR^p,1\right\}
\left[1+\frac{1}{\mu_\mathbf{B}(R)}\right]
\right\}^{-\frac{2}{p}}
\le\nu_\mathbf{B}(R)
\le\mu_\mathbf{B}(R),
\]
which implies that $0<\nu_\mathbf{B}(R)
\le\mu_\mathbf{B}(R)<\infty$. This ends the proof.
\end{proof}

Now, it remains to accomplish the proof of Theorem \ref{thm-1.3}.

\begin{proof}
[\textbf{Proof of Theorem \ref{thm-1.3}}]
Let $u\in\mathcal{D}_{\mathbf{A},p}$ be  satisfying $(ux\cdot\nabla\!_\mathbf{A} u)=|x||\nabla\!_\mathbf{A} u||u|$. It is not difficult to verify that
\begin{align*}
\left(\left|\nabla\!_\mathbf{A} u\right|
+\frac{N-p}{p}\frac{|u|}{|x|}
\right)^{p-2}
|x|^{\frac{2(p-N)}{p}}
\left|\nabla\!_\mathbf{A} \left(u|x|^{\frac{N-p}{p}}\right)\right|^2
=|x|^{p-N}\left|\nabla\!_\mathbf{A} \left(u|x|^{\frac{N-p}{p}}\right)\right|^p,
\end{align*}
then for $1<p<2$,
\begin{align}\label{rc}
&\min\left\{|x|^{p-N}\left|\nabla\!_\mathbf{A} \left(u|x|^{\frac{N-p}{p}}\right)\right|^p,
\left(\frac{N-p}{p}\right)^{p-2}|u|^{p-2}
|x|^{4-p-\frac{2N}{p}}
\left|\nabla\!_\mathbf{A} \left(u|x|^{\frac{N-p}{p}}\right)\right|^2\right\}
\nonumber\\&\quad
=|x|^{p-N}\left|\nabla\!_\mathbf{A} \left(u|x|^{\frac{N-p}{p}}\right)\right|^p.
\end{align}

Based on Lemma \ref{lem-3.3}, let us fix a constant $R>1$ such that $\nu_\mathbf{B}(R)>0$. We define a radially symmetric cut-off function $\chi\in C^{\infty}(\mathbb{R}^N)$ satisfying $0\le\chi\le1$ as follows:
\begin{align*}
\chi=
\begin{cases}
1,
\ \ &\mathrm{if}\ x\in B_R^c(0);\\[0.5mm]
0,
\ \ &\mathrm{if}\ x\in B_{R_2}(0)\setminus B_{R_1}(0);\\[0.5mm]
1,
\ \ &\mathrm{if}\ x\in B_{\frac{1}{R}}(0),
\end{cases}
\end{align*}
where $R_1,R_2$ are two constants satisfying $\frac{1}{R}<R_1<1<R_2<R$. Therefore, $\mathrm{supp}(|\nabla\chi|)
\subset\overline{\mathcal{O}_R}$ and $\mathrm{supp}(1-\chi)
\subset\overline{\mathcal{O}_R}$.

Notice that
\begin{align}\label{3.7}
\int_{\mathbb{R}^N}
&\frac{|u|^p}{|x|^p
\left(1+|\ln|x||^p\right)}\mathrm{d}x
\nonumber\\&=
\int_{\mathbb{R}^N}
\frac{|(1-\chi)u+\chi u|^p}
{|x|^p\left(1+|\ln|x||^p\right)}\mathrm{d}x
\nonumber\\&\le
2^{p-1}\left[\int_{\mathbb{R}^N}
\frac{|(1-\chi)u|^p}
{|x|^p\left(1+|\ln|x||^p\right)}\mathrm{d}x
+\int_{\mathbb{R}^N}
\frac{|\chi u|^p}
{|x|^p\left(1+|\ln|x||^p\right)}\mathrm{d}x
\right]
\nonumber\\&\le
2^{p-1}\bigg[
\int_{\mathcal{O}_R}
\frac{|u|^p}
{|x|^p\left(1+|\ln|x||^p\right)}\mathrm{d}x
+\int_{B_1(0)}
\frac{|\chi u|^p}
{|x|^p\left(1+|\ln|x||^p\right)}\mathrm{d}x
\nonumber\\&\qquad\qquad
+\int_{B_1^c(0)}
\frac{|\chi u|^p}
{|x|^p\left(1+|\ln|x||^p\right)}\mathrm{d}x
\bigg]
\nonumber\\&=
2^{p-1}\left(\Pi_1+\Pi_2+\Pi_3\right),
\end{align}
where
\begin{align*}
\Pi_1&=\int_{\mathcal{O}_R}
\frac{|u|^p}
{|x|^p\left(1+|\ln|x||^p\right)}\mathrm{d}x;\\
\Pi_2&=\int_{B_1(0)}
\frac{|\chi u|^p}
{|x|^p\left(1+|\ln|x||^p\right)}\mathrm{d}x;\\
\Pi_3&=\int_{B_1^c(0)}
\frac{|\chi u|^p}
{|x|^p\left(1+|\ln|x||^p\right)}\mathrm{d}x.
\end{align*}

\textbf{For $\bf\Pi_1$}. We infer from Lemma \ref{lem-3.3} that
{\footnotesize
\begin{align}\label{3.5}
\Pi_1&\le
\int_{\mathcal{O}_R}
\frac{\left|u\right|^p}{|x|^p}
\mathrm{d}x
\nonumber\\&=
\int_{\mathcal{O}_R}
\frac{\left|u|x|^{\frac{N-p}{p}}\right|^p}{|x|^N}
\mathrm{d}x
\nonumber\\&\le
R^N\int_{\mathcal{O}_R}
\left|u|x|^{\frac{N-p}{p}}\right|^p
\mathrm{d}x
\nonumber\\&\le
\frac{R^N}{\nu_\mathbf{B}(R)}\int_{\mathcal{O}_R}
\min\left\{\left|\nabla\!_\mathbf{A}
\left(u|x|^{\frac{N-p}{p}}\right)\right|^p,
\left(\frac{N-p}{p}\right)^{p-2}
|u|^{p-2}|x|^{\frac{(p-2)(N-p)}{p}}
\left|\nabla\!_\mathbf{A}
\left(u|x|^{\frac{N-p}{p}}\right)\right|^2\right\}
\mathrm{d}x
\nonumber\\&=
\frac{R^N}{\nu_\mathbf{B}(R)}\int_{\mathcal{O}_R}
\min\left\{\left|\nabla\!_\mathbf{A}
\left(u|x|^{\frac{N-p}{p}}\right)\right|^p,
|x|^{N-2}\left(\frac{N-p}{p}\right)^{p-2}
|u|^{p-2}|x|^{4-p-\frac{2N}{p}}
\left|\nabla\!_\mathbf{A}
\left(u|x|^{\frac{N-p}{p}}\right)\right|^2\right\}
\mathrm{d}x
\nonumber\\&\le
\frac{\max\left\{R^{2N-p},R^{2N-2}\right\}}
{\nu_\mathbf{B}(R)}
\nonumber\\&\quad\times
\int_{\mathcal{O}_R}
\min\left\{|x|^{p-N}\left|\nabla\!_\mathbf{A} \left(u|x|^{\frac{N-p}{p}}\right)\right|^p,
\left(\frac{N-p}{p}\right)^{p-2}|u|^{p-2}
|x|^{4-p-\frac{2N}{p}}
\left|\nabla\!_\mathbf{A} \left(u|x|^{\frac{N-p}{p}}\right)\right|^2\right\}
\mathrm{d}x
\nonumber\\&=
\frac{R^{2N-p}}{\nu_\mathbf{B}(R)}
\int_{\mathcal{O}_R}
\min\left\{|x|^{p-N}\left|\nabla\!_\mathbf{A} \left(u|x|^{\frac{N-p}{p}}\right)\right|^p,
\left(\frac{N-p}{p}\right)^{p-2}|u|^{p-2}
|x|^{4-p-\frac{2N}{p}}
\left|\nabla\!_\mathbf{A} \left(u|x|^{\frac{N-p}{p}}\right)\right|^2\!\right\}
\mathrm{d}x.
\end{align}}\noindent
As a by-product of \eqref{3.5},
{\footnotesize
\begin{align}\label{3.10}
&\int_{\mathcal{O}_R}
\left|u\right|^p
\mathrm{d}x
\nonumber\\&\ \le
\frac{R^{2N}}
{\nu_\mathbf{B}(R)}
\int_{\mathcal{O}_R}
\min\left\{|x|^{p-N}\left|\nabla\!_\mathbf{A}
\left(u|x|^{\frac{N-p}{p}}\right)\right|^p,
\left(\frac{N-p}{p}\right)^{p-2}
|u|^{p-2}
|x|^{4-p-\frac{2N}{p}}
\!\left|\nabla\!_\mathbf{A}
\left(u|x|^{\frac{N-p}{p}}\right)\right|^2\right\}
\mathrm{d}x.
\end{align}}

\textbf{For $\bf\Pi_2$}. From Lemma \ref{lem-3.1} (with $\widehat{R}=1$), \eqref{rc} and \eqref{1.1}, it follows that
{\footnotesize
\begin{align}\label{3.9}
\Pi_2&\le
\int_{B_1(0)}
\frac{|\chi u|^p}
{|x|^p|\ln|x||^p}\mathrm{d}x
\nonumber\\&\le
\left(\frac{p}{p-1}\right)^p
\int_{B_1(0)}|x|^{p-N}\left|\nabla \left(\chi|u||x|^{\frac{N-p}{p}}\right)\right|^p
\mathrm{d}x
\nonumber\\&=
\left(\frac{p}{p-1}\right)^p
\int_{B_1(0)}
\min\left\{|x|^{p-N}\left|\nabla \left(\chi|u||x|^{\frac{N-p}{p}}\right)\right|^p,
\left(\frac{N-p}{p}\right)^{p-2}
|\chi u|^{p-2}
|x|^{4-p-\frac{2N}{p}}
\left|\nabla \left(\chi|u||x|^{\frac{N-p}{p}}\right)\right|^2\right\}
\mathrm{d}x
\nonumber\\&\le
2\left(\frac{p}{p-1}\right)^p\int_{B_1(0)}
\min\bigg\{\left|\nabla\chi\right|^p|u|^p
+|x|^{p-N}\left|\nabla \left(|u||x|^{\frac{N-p}{p}}\right)\right|^p,
\nonumber\\&\qquad\qquad\qquad\qquad\qquad\quad
\left(\frac{N-p}{p}\right)^{p-2}
|\nabla\chi|^2\frac{|u|^{p}}{|x|^{p-2}}
+\left(\frac{N-p}{p}\right)^{p-2}
|u|^{p-2}
|x|^{4-p-\frac{2N}{p}}\left|\nabla \left(|u||x|^{\frac{N-p}{p}}\right)\right|^2
\bigg\}
\mathrm{d}x
\nonumber\\&\le
2\left(\frac{p}{p-1}\right)^p\int_{B_1(0)}
\max\left\{\left|\nabla\chi\right|^p|u|^p,
\left(\frac{N-p}{p}\right)^{p-2}
\left|\nabla\chi\right|^2\frac{|u|^{p}}{|x|^{p-2}}\right\}
\mathrm{d}x
\nonumber\\&\quad
+2\left(\frac{p}{p-1}\right)^p\int_{B_1(0)}
\min\bigg\{|x|^{p-N}\left|\nabla \left(|u||x|^{\frac{N-p}{p}}\right)\right|^p,
\left(\frac{N-p}{p}\right)^{p-2}
|u|^{p-2}|x|^{4-p-\frac{2N}{p}}\left|\nabla \left(|u||x|^{\frac{N-p}{p}}\right)\right|^2
\bigg\}
\mathrm{d}x
\nonumber\\&\le
2\left(\frac{p}{p-1}\right)^p\int_{B_1(0)}
\left|\nabla\chi\right|^p|u|^p
\mathrm{d}x
+2\left(\frac{p}{p-1}\right)^p
\left(\frac{N-p}{p}\right)^{p-2}
\int_{B_1(0)}
\left|\nabla\chi\right|^2\frac{|u|^{p}}{|x|^{p-2}}
\mathrm{d}x
\nonumber\\&\quad
+2\left(\frac{p}{p-1}\right)^p\int_{B_1(0)}
\min\bigg\{|x|^{p-N}\left|\nabla\!_\mathbf{A} \left(u|x|^{\frac{N-p}{p}}\right)\right|^p,
\left(\frac{N-p}{p}\right)^{p-2}
|u|^{p-2}
|x|^{4-p-\frac{2N}{p}}\left|\nabla\!_\mathbf{A} \left(u|x|^{\frac{N-p}{p}}\right)\right|^2
\bigg\}
\mathrm{d}x.
\end{align}}\noindent
With the aid of \eqref{3.10}, there hold
{\footnotesize
\begin{align*}
&\int_{B_1(0)}\left|\nabla\chi\right|^p|u|^p\mathrm{d}x
\nonumber\\&\quad\le
\|\nabla\chi\|_{L^\infty(\mathbb{R}^N)}^p
\int_{B_1(0)\setminus B_{\frac{1}{R}}(0)}|u|^p\mathrm{d}x
\nonumber\\&\quad\le
\|\nabla\chi\|_{L^\infty(\mathbb{R}^N)}^p
\int_{\mathcal{O}_R}
|u|^p\mathrm{d}x
\nonumber\\&\quad\le
\|\nabla\chi\|_{L^\infty(\mathbb{R}^N)}^p
\frac{R^{2N}}
{\nu_\mathbf{B}(R)}
\int_{\mathcal{O}_R}
\min\left\{|x|^{p-N}\left|\nabla\!_\mathbf{A}
\left(u|x|^{\frac{N-p}{p}}\right)\right|^p,
\left(\frac{N-p}{p}\right)^{p-2}
|u|^{p-2}
|x|^{4-p-\frac{2N}{p}}\left|\nabla\!_\mathbf{A}
\left(u|x|^{\frac{N-p}{p}}\right)\right|^2\right\}
\mathrm{d}x
\end{align*}}\noindent
and
{\footnotesize
\begin{align*}
&\int_{B_1(0)}\left|\nabla\chi\right|^2
\frac{|u|^p}{|x|^{p-2}}\mathrm{d}x
\nonumber\\&\quad\le
\|\nabla\chi\|_{L^\infty(\mathbb{R}^N)}^2
\int_{B_1(0)\setminus B_{\frac{1}{R}}(0)}
\frac{|u|^p}{|x|^{p-2}}\mathrm{d}x
\nonumber\\&\quad\le
\|\nabla\chi\|_{L^\infty(\mathbb{R}^N)}^2
\int_{\mathcal{O}_R}
|u|^p\mathrm{d}x
\nonumber\\&\quad\le
\|\nabla\chi\|_{L^\infty(\mathbb{R}^N)}^2
\frac{R^{2N}}
{\nu_\mathbf{B}(R)}
\int_{\mathcal{O}_R}
\min\left\{|x|^{p-N}\left|\nabla\!_\mathbf{A}
\left(u|x|^{\frac{N-p}{p}}\right)\right|^p,
\left(\frac{N-p}{p}\right)^{p-2}
|u|^{p-2}
|x|^{4-p-\frac{2N}{p}}\left|\nabla\!_\mathbf{A}
\left(u|x|^{\frac{N-p}{p}}\right)\right|^2\right\}
\mathrm{d}x.
\end{align*}}\noindent
Thus, substituting the above two inequalities into \eqref{3.9},
{\footnotesize
\begin{align}\label{3.14}
\Pi_2&\le
2\left(\frac{p}{p-1}\right)^p
\left[\|\nabla\chi\|_{L^\infty(\mathbb{R}^N)}^p
\frac{R^{2N}}
{\nu_\mathbf{B}(R)}
+\left(\frac{N-p}{p}\right)^{p-2}
\|\nabla\chi\|_{L^\infty(\mathbb{R}^N)}^2
\frac{R^{2N}}
{\nu_\mathbf{B}(R)}\right]
\nonumber\\&\qquad
\times\int_{\mathcal{O}_R}
\min\left\{|x|^{p-N}\left|\nabla\!_\mathbf{A}
\left(u|x|^{\frac{N-p}{p}}\right)\right|^p,
\left(\frac{N-p}{p}\right)^{p-2}
|u|^{p-2}
|x|^{4-p-\frac{2N}{p}}\left|\nabla\!_\mathbf{A}
\left(u|x|^{\frac{N-p}{p}}\right)\right|^2\right\}
\mathrm{d}x
\nonumber\\&\quad+
2\left(\frac{p}{p-1}\right)^p\int_{B_1(0)}
\min\left\{|x|^{p-N}\left|\nabla\!_\mathbf{A}
\left(u|x|^{\frac{N-p}{p}}\right)\right|^p,
\left(\frac{N-p}{p}\right)^{p-2}
|u|^{p-2}
|x|^{4-p-\frac{2N}{p}}\left|\nabla\!_\mathbf{A}
\left(u|x|^{\frac{N-p}{p}}\right)\right|^2\right\}
\mathrm{d}x
\nonumber\\&\le
2\left(\frac{p}{p-1}\right)^p
\left[\|\nabla\chi\|_{L^\infty(\mathbb{R}^N)}^p
\frac{R^{2N}}
{\nu_\mathbf{B}(R)}
+\left(\frac{N-p}{p}\right)^{p-2}
\|\nabla\chi\|_{L^\infty(\mathbb{R}^N)}^2
\frac{R^{2N}}
{\nu_\mathbf{B}(R)}+1\right]
\nonumber\\&\qquad
\times\int_{\mathbb{R}^N}
\min\left\{|x|^{p-N}\left|\nabla\!_\mathbf{A}
\left(u|x|^{\frac{N-p}{p}}\right)\right|^p,
\left(\frac{N-p}{p}\right)^{p-2}
|u|^{p-2}
|x|^{4-p-\frac{2N}{p}}\left|\nabla\!_\mathbf{A}
\left(u|x|^{\frac{N-p}{p}}\right)\right|^2\right\}
\mathrm{d}x.
\end{align}}

\textbf{For $\bf\Pi_3$}. Similar arguments to those of $\Pi_2$ show that
{\small
\begin{align}\label{3.11}
\Pi_3&\le
2\left(\frac{p}{p-1}\right)^p
\left[\|\nabla\chi\|_{L^\infty(\mathbb{R}^N)}^p
\frac{R^{2N}}{\nu_\mathbf{B}(R)}
+\left(\frac{N-p}{p}\right)^{p-2}
\|\nabla\chi\|_{L^\infty(\mathbb{R}^N)}^2
\frac{R^{2N+2-p}}
{\nu_\mathbf{B}(R)}+1\right]
\nonumber\\&\quad
\times\int_{\mathbb{R}^N}
\min\left\{|x|^{p-N}\left|\nabla\!_\mathbf{A}
\left(u|x|^{\frac{N-p}{p}}\right)\right|^p,
\left(\!\frac{N-p}{p}\!\right)^{p-2}
|u|^{p-2}
|x|^{4-p-\frac{2N}{p}}\left|\nabla\!_\mathbf{A}
\left(u|x|^{\frac{N-p}{p}}\right)\right|^2\right\}
\mathrm{d}x.
\end{align}}\noindent
Combining \eqref{3.14} and \eqref{3.11}, we get
{\small
\begin{align}\label{3.6}
&\Pi_2+\Pi_3
\nonumber\\&\quad\le
4\left(\frac{p}{p-1}\right)^p
\left[\|\nabla\chi\|_{L^\infty(\mathbb{R}^N)}^p
\frac{R^{2N}}{\nu_\mathbf{B}(R)}
+\left(\frac{N-p}{p}\right)^{p-2}
\|\nabla\chi\|_{L^\infty(\mathbb{R}^N)}^2
\frac{R^{2N+2-p}}
{\nu_\mathbf{B}(R)}+1\right]
\nonumber\\&\qquad
\times\int_{\mathbb{R}^N}
\min\left\{|x|^{p-N}\left|\nabla\!_\mathbf{A}
\left(u|x|^{\frac{N-p}{p}}\right)\right|^p,
\left(\!\frac{N-p}{p}\!\right)^{p-2}
|u|^{p-2}
|x|^{4-p-\frac{2N}{p}}\left|\nabla\!_\mathbf{A}
\left(u|x|^{\frac{N-p}{p}}\right)\right|^2\right\}
\mathrm{d}x.
\end{align}}\noindent
Hence, substituting \eqref{3.5} and \eqref{3.6} into \eqref{3.7}, we can deduce from Theorem \ref{thm-1.2} that
{\small
\begin{align}\label{3.16}
&\int_{\mathbb{R}^N}
\frac{|u|^p}{|x|^p
\left(1+|\ln|x||^p\right)}\mathrm{d}x
\nonumber\\&\quad\le
C_{\mathbf{B},p,N}\int_{\mathbb{R}^N}
\min\bigg\{|x|^{p-N}\left|\nabla\!_\mathbf{A} \left(u|x|^{\frac{N-p}{p}}\right)\right|^p,
\left(\frac{N-p}{p}\right)^{p-2}
|u|^{p-2}
|x|^{4-p-\frac{2N}{p}}
\left|\nabla\!_\mathbf{A} \left(u|x|^{\frac{N-p}{p}}\right)\right|^2\bigg\}
\mathrm{d}x
\nonumber\\&\quad\le
\frac{C_{\mathbf{B},p,N}}{c_3(p)}
\mathcal{H}_{\mathbf{A},p}(u),
\end{align}}\noindent
where
{\footnotesize
\begin{equation}\label{3.17}
C_{\mathbf{B},p,N}
=2^{p-1}\left\{\frac{R^{2N-p}}
{\nu_\mathbf{B}(R)}
+4\left(\frac{p}{p-1}\right)^p
\!\left[\|\nabla\chi\|_{L^\infty(\mathbb{R}^N)}^p
\frac{R^{2N}}{\nu_\mathbf{B}(R)}
+\left(\!\frac{N\!-\!p}{p}\!\right)^{p-2}
\|\nabla\chi\|_{L^\infty(\mathbb{R}^N)}^2
\frac{R^{2N+2-p}}{\nu_\mathbf{B}(R)}
+1\!\right]\!\right\}.
\end{equation}}\noindent
This completes the proof of Theorem \ref{thm-1.3}.
\end{proof}

\begin{rem}
There still exists another version of remainder terms involving logarithmic terms of \eqref{1.6} for $1<p<2$. The detailed proof is stated below.

From \eqref{3.16}, \eqref{1.5} and Theorem \ref{thm-1.1}, if $u\in\mathcal{D}_{\mathbf{A},p}$ satisfies $(ux\cdot\nabla\!_\mathbf{A} u)=|x||\nabla\!_\mathbf{A} u||u|$, we get
\begin{align*}
&\int_{\mathbb{R}^N}
\frac{|u|^p}{|x|^p
\left(1+|\ln|x||^p\right)}\mathrm{d}x
\nonumber\\&\quad\le
C_{\mathbf{B},p,N}\!\int_{\mathbb{R}^N}
\!\min\left\{|x|^{p-N}\left|\nabla\!_\mathbf{A} \!\left(u|x|^{\frac{N-p}{p}}\right)\right|^p,
\!\left(\!\frac{N\!-\!p}{p}\!\right)^{p-2}
|u|^{p-2}
|x|^{4-p-\frac{2N}{p}}
\left|\nabla\!_\mathbf{A} \left(u|x|^{\frac{N-p}{p}}\right)\right|^2\right\}
\mathrm{d}x
\nonumber\\&\quad\le
3^{2-p}C_{\mathbf{B},p,N}\int_{\mathbb{R}^N}
\left(\left|\nabla\!_\mathbf{A} u\right|
+\frac{N-p}{p}\frac{|u|}{|x|}
\right)^{p-2}
|x|^{\frac{2(p-N)}{p}}
\left|\nabla\!_\mathbf{A} \left(u|x|^{\frac{N-p}{p}}\right)\right|^2
\mathrm{d}x
\\&\quad\le
\frac{3^{2-p}C_{\mathbf{B},p,N}}{c_1(p)}
\mathcal{H}_{\mathbf{A},p}(u),
\end{align*}
where $C_{\mathbf{B},p,N}$ and $c_1(p)>0$ are defined by \eqref{3.17} and Theorem \ref{thm-1.1}, respectively.
\end{rem}

\section*{Declarations}

\subsection*{Funding}
\noindent This paper was supported by the National Natural Science Foundation of China (No. 12371120) and Chongqing Graduate Student Research Innovation Project (No. CYB240097).

\subsection*{Data availability statement}
\noindent No data was used for the research described in the article.

\subsection*{Conflict of interest}
\noindent The authors declare no conflict of interest.

\appendix

\section{Appendix: a technical inequality}
\label{sect-appendix}

\noindent In this section, we focus on providing the detailed proof of \eqref{1.4}.

\begin{lem}\label{lem-a}
Let $1<p<2\le N$, for $\eta,\zeta\in\mathbb{C}^N$, we have
\begin{equation*}
\left(|\eta|+|\eta-\zeta|\right)^{p-2}|\zeta|^2
\le\min\left\{|\zeta|^p,|\eta-\zeta|^{p-2}|\zeta|^2\right\}
\le3^{2-p}\left(|\eta|+|\eta-\zeta|\right)^{p-2}|\zeta|^2.
\end{equation*}
\end{lem}

\begin{proof}[\textbf{Proof}]
We continuous the proof in two cases: $|\eta-\zeta|\le|\zeta|$ and $|\eta-\zeta|\ge|\zeta|$.

$\bullet$ If $|\eta-\zeta|\le|\zeta|$, then  $\min\left\{|\zeta|^p,|\eta-\zeta|^{p-2}|\zeta|^2\right\}
=|\zeta|^p$. To complete the proof of this case, it is equivalent to verify that
\begin{equation}\label{2.7}
\left(|\eta|+|\eta-\zeta|\right)^{p-2}
\le|\zeta|^{p-2}
\le3^{2-p}\left(|\eta|+|\eta-\zeta|\right)^{p-2}.
\end{equation}
Notice that
\begin{equation*}
|\eta|=|\eta-\zeta+\zeta|
\le|\eta-\zeta|+|\zeta|
\le2|\zeta|,
\quad
|\eta|=|\eta-\zeta+\zeta|
\ge|\zeta|-|\eta-\zeta|,
\end{equation*}
then
\begin{equation*}
|\eta|+|\eta-\zeta|
\le2|\zeta|+|\eta-\zeta|
\le3|\zeta|,
\ \
|\eta|+|\eta-\zeta|\ge|\zeta|,
\end{equation*}
namely,
\[
\frac{1}{3}\left(|\eta|+|\eta-\zeta|\right)
\le|\zeta|
\le|\eta|+|\eta-\zeta|.
\]
This leads to \eqref{2.7}.

$\bullet$ If $|\eta-\zeta|\ge|\zeta|$, then $\min\left\{|\zeta|^p,|\eta-\zeta|^{p-2}|\zeta|^2\right\}
=|\eta-\zeta|^{p-2}|\zeta|^2$. Given this, it remains to confirm that
\begin{equation}\label{2.8}
\left(|\eta|+|\eta-\zeta|\right)^{p-2}
\le|\eta-\zeta|^{p-2}
\le3^{2-p}\left(|\eta|+|\eta-\zeta|\right)^{p-2}.
\end{equation}
Observe that
\begin{equation*}
|\eta|=|\eta-\zeta+\zeta|
\le|\eta-\zeta|+|\zeta|
\le2|\eta-\zeta|,
\end{equation*}
then
\begin{equation*}
|\eta|+|\eta-\zeta|\le3|\eta-\zeta|,
\ \
|\eta-\zeta|\le|\eta|+|\eta-\zeta|,
\end{equation*}
that is,
\[
\frac{1}{3}\left(|\eta|+|\eta-\zeta|\right)
\le|\eta-\zeta|
\le|\eta|+|\eta-\zeta|.
\]
This gives \eqref{2.8}. The proof is completed.
\end{proof}

\end{document}